\documentclass[pdflatex,sn-mathphys-num, iicol]{sn-jnl}

\usepackage{graphicx}%
\usepackage{multirow}%
\usepackage{amsmath,amssymb,amsfonts}%
\usepackage{amsthm}%
\usepackage{mathrsfs}%
\usepackage[title]{appendix}%
\usepackage{xcolor}%
\usepackage{textcomp}%
\usepackage{manyfoot}%
\usepackage{booktabs}%
\usepackage{listings}%
\usepackage{empheq}
\usepackage{mathtools}
\usepackage{stfloats}
\usepackage[ruled,vlined]{algorithm2e}
\usepackage{subcaption}

\theoremstyle{thmstyleone}%
\theoremstyle{thmstyletwo}%
\newtheorem{remark}{Remark}%

\theoremstyle{thmstylethree}%

\newcommand{\blockDist}{{\textsc{BlockDist}}}
\newcommand{\patAlloc}{{\textsc{PatAlloc}}}
\newcommand{\block}{{\textsc{BlockDist-PatAlloc}}}
\newcommand{\milp}{{\textsc{MILP-PatAlloc}}}

\begin{document}

\title[Two-Stage OR Allocation Framework]{A Two-Stage Operating Room Allocation Framework for Reducing Surgical 
Waiting Lists in Public Hospitals}


\author*[1]{\fnm{Jhoan} \sur{Báez}}\email{sbaez@cmm.uchile.cl}

\author[2]{\fnm{Juan} \sur{Hoyos}}\email{jhoyoss@unal.edu.co}

\author[1]{\fnm{Víctor} \sur{Riquelme}}\email{vriquelme@dim.uchile.cl}



\author[3]{\fnm{Michel} \sur{Royer}}\email{mroyer@calvomackenna.cl}

\author[1,4]{\fnm{Héctor} \sur{Ramírez}}\email{hramirez@dim.uchile.cl}

\affil*[1]{\orgdiv{Center for Mathematical Modeling}, \orgname{Universidad de Chile and CNRS (IRL 2807)}, \orgaddress{ \city{Santiago},  \country{Chile}}}

\affil[2]{\orgname{Universidad Nacional de Colombia}, \orgaddress{ \city{La Paz, Cesar},  \country{Colombia}}}


\affil[3]{ \orgname{Dr. Luis Calvo Mackenna Hospital}, \orgaddress{ \city{Santiago},  \country{Chile}}}

\affil[4]{\orgdiv{Department of Mathematical Engineering}, \orgname{Universidad de Chile}, \orgaddress{ \city{Santiago},  \country{Chile}}}


\abstract{
Operating-room allocation is a major challenge for public hospitals with limited surgical capacity and large elective waiting lists.  This work proposes a two-stage operating-room allocation framework for specialty-block scheduling environments. 
The methodology combines  a mixed-integer linear programming model for medical specialty block allocation with a priority-based patient allocation procedure. The  framework was evaluated through one-week and multi-week simulation scenarios using parameters estimated from historical surgical and  waiting-list data. The proposed methodology was compared against an integrated mixed-integer linear programming baseline adapted from  the literature under realistic operational disruptions, including failed patient confirmations and surgery suspensions. 
Results show that  the proposed framework achieves a more balanced distribution between offered and demanded surgical time across specialties while maintaining  high operating-room utilization and substantially lower computational times than the integrated baseline. 
In the 20-week scenario, the  proposed methodology achieved lower waiting times and operated a larger number of patients while preserving operational flexibility. These  results suggest that decomposition-based operating-room allocation strategies provide a practical and computationally efficient alternative  for reducing surgical waiting lists in public hospitals.
}





\maketitle

\section{Introduction}\label{introduction}

Operating rooms are among the most resource-intensive and capacity-constrained components of modern  healthcare systems. Their efficient allocation directly impacts patient access to care, hospital  productivity, healthcare costs, and quality of service. Consequently, operating room planning and  scheduling has been extensively studied in operations research and healthcare optimization through  mathematical programming, stochastic optimization, simulation, heuristics, and hybrid approaches  \citep{cardoen2010operating,guerriero2011operational,zhu2019operating,samudra2016}.


The literature on operating-room scheduling has primarily focused on objectives such as operating-room utilization, overtime minimization, throughput maximization, and operational cost reduction. Public healthcare systems, however, often operate under a different set of priorities, where reducing surgical waiting lists under severe resource limitations becomes a central objective. In highly congested environments, long waiting times may directly affect patient outcomes, quality of care, and healthcare equity. Consequently, reducing the mean and median waiting times of elective surgery patients represents a critical operational challenge.

This problem is particularly relevant in the Chilean public healthcare system, where elective  surgical waiting lists constitute a persistent national challenge. Public hospitals must  continuously allocate limited operating-room capacity among multiple competing medical  specialties while handling heterogeneous surgical durations, restricted surgical-team  availability, evolving waiting lists, and frequent operational disruptions. In practice, many  scheduling decisions are still performed manually using historical allocation patterns and  expert knowledge, which may generate important imbalances between the actual surgical demand  and the offered operating-room capacity across specialties.


Integrated operating-room scheduling models have shown promising theoretical results and constitute an important direction in the operating-room scheduling literature (see \cite{al2025comprehensive}). Their applicability in public hospitals, however, is often influenced by the scale of the waiting list and the dynamic nature of daily operations. As the number of patients increases, integrated formulations may become computationally demanding, while the resulting schedules can be difficult to update in response to operational changes.

Public hospitals typically operate in highly dynamic environments, where elective surgery schedules are continuously influenced by patient cancellations, unsuccessful confirmations, emergency reallocations, missing examinations, and last-minute modifications. Under these conditions, scheduling methodologies must not only produce high-quality solutions but also allow rapid recomputation and straightforward schedule adjustments. These operational requirements motivate the development of flexible decomposition-based approaches such as the one proposed in this work.

The proposed framework was developed collaboratively with Hospital Dr. Luis Calvo Mackenna (HLCM), a high-complexity pediatric public hospital in Santiago, Chile, using anonymized historical  elective-surgery and waiting-list data. Unlike fully integrated operating-room scheduling  approaches, the proposed methodology explicitly prioritizes operational flexibility, fast  execution times, and adaptability to real hospital workflows. In particular, the framework  preserves the specialty-block structure required by the hospital while allowing efficient  patient reassignment and schedule modifications under operational disruptions.

To evaluate the proposed methodology, we compare the proposed framework against an integrated  mixed-integer linear programming baseline inspired by \cite{wolff2012}, adapted to prioritize  waiting-list reduction instead of overtime minimization. Computational experiments based on  real hospital data show that the proposed methodology achieves substantial reductions in mean  and median waiting times while maintaining practical utilization rates and significantly lower  computational times than the integrated optimization baseline.

The main contributions of this work are summarized as follows:
\begin{itemize}
\item We propose a practical two-stage operating-room allocation framework specifically designed for public hospitals operating under specialty-block scheduling policies.
\item We introduce a block-distribution methodology that balances offered and demanded operating-room time across specialties while preserving operational flexibility.
\item We incorporate practical scheduling disruptions, including failed patient confirmations and surgery suspensions, into the evaluation framework.
\item We evaluate the proposed methodology using real data from a high-complexity pediatric public hospital in Chile.
\item We demonstrate that the proposed framework achieves substantial reductions in surgical waiting times while requiring significantly lower computational times than an integrated MILP-based approach.
\end{itemize}

The remainder of this paper is organized as follows. Section~\ref{sec:literature} reviews the  related literature on operating-room scheduling and surgical waiting-list management.  Section~\ref{sec:context} describes the hospital setting and the operational characteristics  motivating the proposed approach. Section~\ref{sec:framework} presents the proposed two-stage  operating-room allocation framework and the associated optimization models.  Section~\ref{sec:experimental_design} describes the experimental setup, simulation  methodology, and evaluation metrics. Section~\ref{sec:results} presents the computational results. Finally, Sections~\ref{sec:discussion} and~\ref{sec:conclusion} provide the discussion and conclusions of the study.

\section{Hospital Context and Problem Description}\label{sec:context}

HLCM performs elective surgeries using a specialty-block scheduling policy, where operating-room time is divided into morning (AM) and afternoon (PM) blocks and each block is assigned to a single medical specialty. Although this structure simplifies the coordination of surgeons, anesthesiologists, nurses, and surgical equipment, it also introduces operational constraints that are not always considered in classical operating-room scheduling formulations.

The hospital has eight operating rooms dedicated to elective surgeries, one of which is reserved exclusively for cardiovascular procedures. Elective surgeries are scheduled weekly from Monday to Friday, resulting in a weekly capacity of 39 hours per operating room. Due to limitations in surgeons and anesthesiologists, it is uncommon for the hospital to operate more than two simultaneous blocks of the same specialty.

Historically, operating-room time allocation at HLCM was performed manually using historical scheduling patterns and operational experience. Under this approach, specialties tend to receive relatively fixed proportions of weekly operating-room time independently of the actual demand observed in the waiting list. Consequently, some specialties accumulate disproportionately large waiting lists and prolonged waiting times, while others receive more operating-room time than required.

The database used in this study contains anonymized historical information from elective surgeries  and waiting lists collected over 26 consecutive months. To isolate the impact of the scheduling decisions, surgery durations are treated as deterministic and estimated using the historical average duration of interventions in which the corresponding procedure appears as the principal surgery. Since multiple procedures may be performed during a single intervention, the principal procedure provides a consistent reference for duration estimation across all surgery types.
In addition to the surgical duration, a cleaning and preparation time of $u=15$ minutes is considered between consecutive surgeries. This value was not estimated from the historical records but was established by hospital staff based on their operational experience and standard clinical practice.

A distinctive feature of the hospital's scheduling policy is that operating-room capacity is allocated at the specialty level rather than at the individual surgeon level. Consequently, the proposed methodology focuses on assigning operating-room blocks to specialties, leaving the internal assignment of surgeons to the hospital. This design choice was explicitly requested by hospital management, as operating-room allocation is considered an institutional planning problem rather than an evaluation of individual physician performance. Although previous studies have shown that surgeon-specific characteristics can improve surgery-duration predictions, the hospital prioritizes a specialty-based planning strategy that promotes operational robustness, facilitates schedule modifications, and avoids decisions directly linked to individual practitioners.

The waiting list analyzed in this work presents substantial heterogeneity among specialties in terms of waiting-list size and waiting-time distributions. Moreover, elective scheduling is frequently affected by cancellations, failed patient confirmations, and operational disruptions, making flexibility and rapid schedule adaptation essential characteristics for any practically deployable scheduling methodology.

Motivated by these operational requirements, this work proposes a two-stage operating-room allocation framework specifically designed for public hospitals operating under specialty-block scheduling policies and continuously evolving waiting lists.


\section{Literature Review}\label{sec:literature}

Operating room planning and scheduling has been extensively studied in operations research due to its impact on healthcare efficiency, resource utilization, patient access, and hospital costs \cite{cardoen2010operating,guerriero2011operational,zhu2019operating}. Existing approaches cover strategic, tactical, and operational decision levels, and include methodologies based on mixed-integer programming, stochastic optimization, simulation, heuristics, and robust optimization \cite{samudra2016,Rahimi2020}.

A substantial portion of the literature focuses on integrated optimization models that simultaneously determine operating-room allocation, patient assignment, and surgery sequencing. Several Mixed-Integer Linear Programming (MILP) formulations have been proposed to maximize utilization and minimize overtime while considering resource constraints and downstream dependencies \cite{Santibanez2007,Latorre2016,ThomasSchneider2020}. Other works incorporate balancing mechanisms between offered and demanded operating-room time \cite{Roshanaei2020b,Zhou2020}, as well as stochastic and robust formulations addressing uncertain surgery durations and emergency arrivals \cite{Kamran2018,Wang2020,Atighehchian2020}.


Fully integrated scheduling models provide a powerful optimization framework for jointly addressing multiple operating-room scheduling decisions. Their implementation in real hospital environments, however, is often influenced by computational requirements and the need to accommodate frequent operational changes. In practice, elective surgery schedules are continuously affected by cancellations, failed patient confirmations, emergency reallocations, and staff availability changes. These characteristics have motivated the development of decomposition-based methodologies that separate tactical block allocation from operational patient scheduling \cite{Jebali2006,Luo2019,Aringhieri2015}. Complementary approaches based on heuristics, metaheuristics, decomposition techniques, and reinforcement learning have also been proposed to improve scalability and adaptability \cite{Rizk2011,Roshanaei2021,dolatkhah2026reinforcement}.

Another important research direction concerns uncertainty in surgery durations. Inaccurate duration estimates may generate overtime or idle operating-room time, motivating the development of stochastic scheduling models and predictive methodologies \cite{denton2007,master2017,bartek2018,edelman2017}. However, accurate surgery-duration estimation remains a challenging problem, particularly in public healthcare systems with limited operational information.

Within the Chilean public healthcare context, reducing elective surgical waiting lists has become a major operational priority. In particular, \cite{wolff2012} proposed integer programming formulations for HLCM focused on maximizing operating-room occupation while penalizing overtime. Additional Chilean contributions include stochastic and chance-constrained scheduling approaches for public hospitals \cite{Azar2017,Azar2022,barrera2020}.

Most existing studies prioritize objectives such as utilization maximization, overtime reduction, or throughput optimization. In contrast, this work focuses primarily on reducing mean and median waiting times under the operational constraints of a public hospital operating with specialty-based block scheduling policies. 

\section{Two-Stage Operating Room Allocation Framework}
\label{sec:framework}

In this section, we present the proposed two-stage operating room allocation framework designed to support surgical scheduling decisions in public hospitals under persistent waiting-list pressure. The framework aims to reduce mean and median patient waiting times while respecting operational constraints associated with operating-room availability and specialty-specific surgical resources.

The proposed methodology decomposes the scheduling process into two sequential stages operating at different decision levels. The first stage corresponds to a tactical operating-room block allocation problem, where available operating-room blocks are assigned to medical specialties over a finite planning horizon. The second stage corresponds to an operational patient scheduling procedure, where patients are assigned to the previously allocated blocks according to hospital-defined prioritization criteria.

This decomposition provides several practical advantages. First, it separates long-term capacity balancing decisions from short-term patient scheduling decisions, reducing the computational complexity of the overall problem. Second, it increases operational flexibility, since modifications in patient assignments can be performed without recomputing the complete operating-room allocation. Finally, the framework is consistent with real hospital workflows, where operating-room capacity allocation and patient confirmation are typically managed independently.

The proposed framework consists of the following three components:

\begin{enumerate}
    \item \textbf{Patient prioritization:} the waiting list is partitioned according to medical specialties, and patients are ranked according to a hospital-defined priority criterion. This prioritization may depend on clinical urgency, accumulated waiting time, diagnosis,  or combinations thereof.
    
    \item \textbf{Operating-room block allocation:} available operating-room blocks are assigned to medical specialties while balancing the relationship between demanded surgical time and     offered operating-room capacity under operational constraints.
    
    \item \textbf{Patient scheduling:} patients are assigned to the blocks corresponding to their specialty while respecting the previously defined prioritization order and operating-room time constraints.
\end{enumerate}

The present work focuses on the second and third stages of the framework. In particular, we  propose a mixed integer linear programming model for the operating-room block allocation problem and a priority-based heuristic for patient scheduling.

\subsection{General notation}

Let $M$ denote the set of medical specialties and $P$ the set of patients in the waiting list. For each specialty $m\in M$, let $P_m\subseteq P$ denote the subset of patients associated with specialty $m$, such that
\[
P = \bigcup_{m\in M} P_m.
\]

Each patient is associated with a principal surgical specialty, determined during the prioritization process. This assignment allows each patient to be uniquely considered within the scheduling procedure.

Let $R$ denote the set of available operating rooms and $T$ the set of days within the planning horizon. In practice, the planning horizon typically corresponds to one business week, although the framework can be adapted to other scheduling horizons.

For each patient $p\in P$, let:
\begin{itemize}
    \item $d_p$ denote the estimated duration of the surgery associated with patient $p$;
    \item $w_p$ denote the current waiting time of patient $p$ in the waiting list;
    \item $k_p$ denote the priority score assigned to patient $p$.
\end{itemize}

The estimated surgical duration $d_p$ is obtained from historical average procedure durations associated with the corresponding surgical intervention. Although more sophisticated predictive models for surgery duration estimation may be incorporated in future implementations, the present work focuses exclusively on the operating-room allocation and patient scheduling components. Importantly, the benchmark methodology used for comparison relies on the same duration estimation assumption, ensuring a fair comparison between approaches.

Each day $t\in T$ is divided into two operating shifts:
\begin{itemize}
    \item a morning shift (AM) of duration $q_t^{\rm AM}$ hours;
    \item an afternoon shift (PM) of duration $q_t^{\rm PM}$ hours.
\end{itemize}

An operating-room block is defined as the tuple
\[
(t,r,\text{shift}),
\]
where $t\in T$ corresponds to the operating day, $r\in R$ corresponds to the operating room, 
and $\text{shift}\in\{\text{AM},\text{PM}\}$ corresponds to the operating shift.

The notation used throughout the proposed framework and the benchmark formulation is 
summarized in Table~\ref{tabla_constantes}. Since both the \blockDist\ and \milp\ methodologies 
share most of the operational and scheduling variables, we introduce the notation here to 
avoid redundancy and to improve the readability of the subsequent formulations.

\begin{table}[htbp]
\centering
\caption{Notation and description of model parameters used in the proposed operating-room 
allocation framework.}
\label{tabla_constantes}

\footnotesize
\setlength{\tabcolsep}{3pt}
\renewcommand{\arraystretch}{1.15}

\begin{tabular}{p{0.18\columnwidth}
                p{0.48\columnwidth}
                p{0.22\columnwidth}}
\toprule
Parameter & Description & Range \\
\midrule
$M$ & Set of medical specialties & Finite set \\[0.5mm]
$m$ & Medical specialty & $m\in M$ \\[0.5mm]
$P$ & Set of patients in waiting list & Finite set \\[0.5mm]
$P_m$ & Set of patients in waiting list of specialty $m$ & $P_m\subset P$ \\[0.5mm]
$p$ & Patient in waiting list & $p\in P$ \\[0.5mm]
$R$ & Set of operating rooms & Finite set \\[0.5mm]
$r$ & Operating room & $r\in R$ \\[0.5mm]
$T$ & Set of Planning horizon days & Finite set \\[0.5mm]
$t$ & Day index in the scheduling & $t\in T$ \\[0.5mm]
$q_t^{\rm AM}$ & Time length of morning shift on day $t$ (in hours) & $q_t^{\rm AM}>0$ \\[0.5mm]
$q_t^{\rm PM}$ & Time length of afternoon shift on day $t$ (in hours) & $q_t^{\rm PM}>0$ \\[0.5mm]
$d_p$ & Duration of surgery of patient $p$ (in minutes) & $d_p>0$ \\[0.5mm]
$k_p$ & Priority of patient $p$ (dimensionless) & $k_p>0$ \\[0.5mm]
$w_p$ & Waiting time of patient $p$ in waiting list (in days) & $w_p\geq0$ \\[0.5mm]
$u$ & Cleaning time of room after surgery (in minutes) & $u>0$ \\[0.5mm]
$N_{mt}^{\rm AM}$ & Maximum simultaneous rooms for specialty $m$ on shift AM of day $t$ & $N_{mt}^{\rm AM}\in\mathbb{Z}_{\geq 0}$ \\[0.5mm]
$N_{mt}^{\rm PM}$ & Maximum simultaneous rooms for specialty $m$ on shift PM of day $t$ & $N_{mt}^{\rm PM}\in\mathbb{Z}_{\geq 0}$ \\
$b_{mt}^{\rm AM}$ & Number of morning blocks on day $t$ assigned to specialty $m$ &  $b_{mt}^{\rm AM}\in\mathbb{Z}_{\geq 0}$  \\[0.5mm]
$b_{mt}^{\rm PM}$ & Number of afternoon blocks on day $t$ assigned to specialty $m$ &  $b_{mt}^{\rm PM}\in\mathbb{Z}_{\geq 0}$  \\
\bottomrule
\end{tabular}
\end{table}

\subsection{Operating-room block allocation model}

The first stage of the proposed framework corresponds to the tactical allocation of operating-room blocks to medical specialties. This stage is performed through the \blockDist\ methodology, which internally solves a mixed integer linear programming problem to determine how the available operating-room capacity should be distributed among specialties during the planning horizon.

The main objective of this stage is to balance the relationship between:
\begin{itemize}
    \item the surgical demand associated with each specialty, measured through the estimated surgical time required by the waiting list;
    \item and the operating-room time effectively assigned to that specialty.
\end{itemize}

This balancing mechanism aims to reduce the accumulation of patients in underserved specialties while avoiding inefficient over-allocation of operating-room capacity. 
Unlike formulations focused exclusively on utilization or overtime minimization, the proposed approach explicitly incorporates backlog pressure into the allocation process.

The allocation process must also satisfy operational constraints related to:
\begin{itemize}
    \item operating-room availability;
    \item simultaneous specialty capacity limits;
    \item surgeon availability;
    \item anesthesiology and nursing constraints.
\end{itemize}

For each specialty $m\in M$ and day $t\in T$, let
\[
b_{mt}^{\rm AM}, \quad b_{mt}^{\rm PM}
\]
denote the number of operating-room blocks assigned to specialty $m$ during the morning and afternoon shifts, respectively.

Since operating rooms are operationally interchangeable during the tactical allocation stage, the decision variables represent the number of blocks assigned to each specialty rather than explicit room assignments.

\bmhead{Demand and capacity indicators}

To define the allocation model, we first introduce several quantities associated with the available operating-room capacity and the surgical demand observed in the waiting list.

The total operating-room capacity available during the planning horizon is defined as
\[
H := \sum_{t\in T} |R| \left(q_t^{\rm AM}+q_t^{\rm PM}\right).
\]

For each specialty $m\in M$, the total demanded surgical time associated with the waiting list is defined as
\[
D_m := \sum_{p\in P_m} d_p,
\]
while the total demanded surgical time across all specialties is
\[
D := \sum_{p\in P} d_p.
\]

Using these quantities, we define the ideal proportional operating-room capacity associated with specialty $m$ as
\[
h_m := H \cdot \frac{D_m}{D},
\]
which corresponds to the fraction of the total available operating-room time that specialty $m$ should receive if the operating-room capacity were distributed proportionally to the demanded surgical time in the waiting list.

The total operating-room time assigned to specialty $m$ during the planning horizon is given by
\[
Z_m :=
\sum_{t\in T}
q_t^{\rm AM} b_{mt}^{\rm AM}
+
q_t^{\rm PM} b_{mt}^{\rm PM}.
\]


\bmhead{Objective function}

The objective of the block-allocation stage is to distribute the available operating-room capacity among medical specialties according to their relative priority. Since patient assignment is performed later during the patient-allocation stage, this model determines how much operating-room time should be allocated to each specialty during the planning horizon.

To represent different prioritization policies, each specialty $m\in M$ is assigned a nonnegative weight $\omega_m$, reflecting its relative importance in the allocation process. The model maximizes the weighted operating-room time assigned to specialties:
\begin{equation}
\label{eq:blockdist_obj}
\max_{(b_{mt}^{\rm AM},b_{mt}^{\rm PM})}
\quad
\sum_{m\in M}
\sum_{t\in T}
\omega_m
\left(
q_t^{\rm AM} b_{mt}^{\rm AM}
+
q_t^{\rm PM} b_{mt}^{\rm PM}
\right).
\end{equation}

Since the total operating-room capacity is fixed, the objective function does not increase the overall amount of available surgical time. Instead, it determines how that capacity is distributed among specialties according to the selected prioritization criterion.

In this work, specialty priorities are determined using the average waiting time of patients associated with each specialty:
\begin{equation}
\omega_m =
\frac{1}{|P_m|}
\sum_{p\in P_m} w_p,
\quad m\in M.
\end{equation}

Under this definition, specialties whose patients experience longer waiting times receive higher priority in the allocation process. Consequently, the model tends to assign more operating-room blocks to specialties with greater waiting-list pressure. Moreover, because the arithmetic mean gives greater influence to exceptionally long waiting times, specialties containing highly delayed patients naturally receive additional priority.

More generally, the framework allows alternative definitions of specialty weights according to hospital objectives. For example, weights may be defined using the proportional demanded surgical time,
\begin{equation}
\omega_m = h_m,
\quad m\in M,
\end{equation}
or through hospital-defined urgency or clinical prioritization policies.

Therefore, the proposed framework can support different operational strategies, such as prioritizing specialties with large waiting lists, responding to urgent clinical needs, or allocating capacity proportionally to observed demand.

\begin{remark}
In the numerical experiments presented in this work, patient prioritization is based on waiting time in the waiting list (FIFO criterion), due to the absence of historical urgency-score records. Nevertheless, the proposed framework naturally allows the incorporation of hospital-defined clinical priority scores.
\end{remark}

\bmhead{Constraints}

The allocation model is subject to the following operational constraints.

\begin{enumerate}

\item \textbf{Specialty capacity constraints.}

The number of simultaneously assigned operating rooms for specialty $m$ cannot exceed the specialty-specific operational capacity during each shift:
\begin{equation}
0
\leq
b_{mt}^{\rm AM}
\leq
N_{mt}^{\rm AM},
\qquad
\forall m\in M,\ t\in T,
\end{equation}
\begin{equation}
0
\leq
b_{mt}^{\rm PM}
\leq
N_{mt}^{\rm PM},
\qquad
\forall m\in M,\ t\in T.
\end{equation}

These constraints represent practical limitations associated with surgeon availability, anesthesiology teams, nursing staff, and specialty-specific operational restrictions.

\item \textbf{Operating-room assignment constraints.}

Every available operating-room block must be assigned to exactly one specialty:
\begin{equation}
\sum_{m\in M}
b_{mt}^{\rm AM}
=
|R|,
\qquad
\forall t\in T,
\end{equation}
\begin{equation}
\sum_{m\in M}
b_{mt}^{\rm PM}
=
|R|,
\qquad
\forall t\in T.
\end{equation}

Since each block is associated with a unique specialty, mixed-specialty blocks are not allowed.

\item \textbf{Demand-capacity balancing constraints.}

The total operating-room time assigned to each specialty should remain close to its ideal proportional demand $h_m$. To enforce this balance, we impose the constraints
\begin{equation}
\label{eq:balance_constraint}
z_{\rm low}(h_m)
\leq
Z_m
\leq
z_{\rm up}(h_m),
\qquad
\forall m\in M,
\end{equation}
where:
\begin{itemize}
    \item $z_{\rm low}(h_m)$ denotes the largest feasible lower bound compatible with the block structure;
    
    \item $z_{\rm up}(h_m)$ denotes the smallest feasible upper bound compatible with the block structure.
\end{itemize}

These bounds account for the fact that operating-room assignments are performed using indivisible AM and PM blocks, making exact proportional allocations generally infeasible.

\end{enumerate}

\bmhead{\blockDist\ formulation}

The mixed integer linear programming problem associated with the \blockDist\ methodology is presented in Problem~(\blockDist).

\begin{figure*}
\begin{empheq}[left={{{(\blockDist)}}\quad\empheqlbrace\quad}]{align*}
\quad\max \quad&  \sum_{m\in M \atop t \in T} \omega_m (q_{t}^{{\rm AM}} b_{mt}^{{\rm AM}} + q_{t}^{{\rm PM}} b_{mt}^{{\rm PM}}),  \\[1mm]
s.t. \quad& \sum_{m\in M} b_{mt}^{{\rm AM}} = |R| \,\quad \forall t\in T, \\[1mm]
&\sum_{m\in M} b_{mt}^{{\rm PM}} = |R| \,\quad \forall t\in T,\\[1mm]
&  0\leq b_{mt}^{{\rm AM}}\leq N_{mt}^{{\rm AM}} \,\quad \forall t\in T,\quad \forall m\in M,\\[1mm]
& 0\leq b_{mt}^{{\rm PM}}\leq N_{mt}^{{\rm PM}} \,\quad \forall t\in T,\quad \forall m\in M,\\[1mm]
& z_{low}(h_m) \,\leq\, Z_m \,\leq\, z_{up}(h_m),\quad \forall m\in M.
\end{empheq}
\end{figure*}

\begin{remark}

A natural alternative formulation would consist in directly minimizing the absolute difference between the offered and demanded operating-room time for each specialty:
\[
\min \sum_{m\in M} |h_m - Z_m|.
\]

However, such formulation may lead to undesirable operational behavior in specialties with relatively small waiting lists but highly urgent patients. In particular, specialties with low aggregate demand could receive no operating-room blocks, even when some patients require urgent surgical intervention or have experienced excessively long waiting times.

The proposed weighting mechanism avoids this issue by allowing hospitals to explicitly incorporate urgency or prioritization criteria into the operating-room allocation process.

\end{remark}

\subsection{Patient scheduling procedure}

Once the operating-room blocks have been assigned to specialties, the second stage of the framework corresponds to the operational scheduling of patients into the allocated blocks. This stage is performed through the \patAlloc\ procedure.

At this stage, each operating-room block is associated with exactly one medical specialty. The objective of the patient scheduling procedure is therefore to assign patients from the corresponding specialty waiting list into the available operating-room blocks while respecting hospital prioritization criteria and operational time constraints.

Patients associated with specialty $m$ are first sorted according to the selected prioritization criterion. Although the framework supports arbitrary hospital-defined prioritization scores, the experiments presented in this work consider waiting-time prioritization (FIFO scheduling), due to the absence of historical urgency-score records.

The patient scheduling procedure follows a greedy priority-based assignment strategy. Patients are sequentially assigned to the available blocks associated with their specialty until the remaining capacity of the block is exhausted.

Three operational constraints must be satisfied during this process:

\begin{enumerate}

\item A cleaning time $u$ must be considered after each surgery.

\item The total surgical and cleaning time assigned to a block cannot exceed the block duration.

\item A patient cannot be assigned to more than one block during the planning horizon.

\end{enumerate}

Let
\[
x_{prt}^{\rm shift}
\]
denote a binary variable indicating whether patient $p$ is assigned to operating room $r$ on day $t$ during the corresponding shift.

We additionally define
\[
O_p =
\sum_{t\in T}
\sum_{r\in R}
\left(
x_{prt}^{\rm AM}
+
x_{prt}^{\rm PM}
\right),
\]
which indicates whether patient $p$ has already been assigned during the planning horizon.

The patient scheduling procedure is summarized in Algorithm~\ref{alg:patient_allocation}.

\begin{algorithm}[htbp]
\caption{Patient scheduling procedure (\patAlloc)}
\label{alg:patient_allocation}

\small

\KwResult{Assignment of patients to operating-room blocks}

Initialize
\[
(x_{prt}^{\rm shift})_{(p,r,t,\rm shift)} \gets 0
\]
and
\[
(O_p)_{p\in P} \gets 0.
\]

\For{$t \in T$}{
    
    \For{$\rm shift \in \{AM,PM\}$}{
        
        \For{$r \in R$}{
            
            $m \gets$ specialty assigned to block $(t,r,\rm shift)$\;
            
            $c \gets q_t^{\rm shift}$\;
            
            Sort patients in $P_m$ according to decreasing priority\;
            
            \For{$p \in P_m$}{
                
                \If{$c \geq d_p + u$ \textbf{and} $O_p = 0$}{
                    
                    $x_{prt}^{\rm shift} \gets 1$\;
                    
                    $O_p \gets 1$\;
                    
                    $c \gets c - (d_p + u)$\;
                }
            }
        }
    }
}

\end{algorithm}

The proposed two-stage framework provides important operational advantages in practice. Since the tactical allocation stage is separated from the patient scheduling stage, modifications in patient assignments can be performed without recomputing the complete operating-room allocation.

This flexibility is particularly relevant in public hospital environments, where patient confirmations and cancellations occur frequently. In practice, hospital personnel must contact patients after scheduling to confirm the intervention date, and unsuccessful confirmations are common. Consequently, patient replacement mechanisms must remain operationally simple and flexible.

Moreover, since operating rooms assigned to specialties during the tactical stage are operationally interchangeable under compatible capacity constraints, specialty blocks can be rearranged when convenient for hospital personnel. For example, blocks assigned to the same specialty on different days may be grouped together to improve operational continuity for surgical teams.

Finally, when waiting-time prioritization is adopted, the proposed scheduling procedure naturally promotes reductions in both the mean and median waiting times of the remaining waiting list, particularly when patients with long waiting times are systematically prioritized for assignment.

\section{Experimental Design}\label{sec:experimental_design}

In this section, we describe the experimental setting used to evaluate the proposed two-stage operating-room allocation framework (\block). The experiments were conducted using historical information from Dr. Luis Calvo Mackenna Hospital (HLCM), a public pediatric hospital in Santiago, Chile. The proposed framework was compared against a benchmark patient-level optimization approach, denoted by \milp, described in Appendix~\ref{MILP}. 

The objective of the experimental analysis is to evaluate the ability of the proposed framework to:
\begin{itemize}
    \item reduce mean and median waiting times;
    \item increase the number of operated patients;
    \item maintain high operating-room utilization;
    \item and improve the balance between offered and demanded surgical time across specialties.
\end{itemize}

\subsection{Hospital setting and operational assumptions}

HLCM has eight operating rooms dedicated to elective surgeries, one of which is reserved exclusively for cardiovascular interventions. The remaining operating rooms are shared among the other surgical specialties according to operational priorities and resource availability.

Elective surgeries are scheduled during weekdays, from Monday to Friday. The planning horizon therefore consists of five working days: $|T| = 5$.

From Monday to Thursday, the daily schedule is divided into: one morning (AM) shift of 5 hours (08:00--13:00), and one afternoon (PM) shift of 3 hours (14:00--17:00).
On Fridays, the PM shift is reduced to 2 hours (14:00--16:00). Consequently, each operating room provides a weekly capacity of 39 hours.

A cleaning time of
\[
u = 15 \text{ minutes}
\]
is considered after each surgery.

Overtime is not allowed during the scheduling process, except in cases where the estimated duration of a surgery exceeds the capacity of the largest available block (5 hours). In such situations, both AM and PM blocks are assigned consecutively to the corresponding specialty.

Due to limitations in surgeon and anesthesiology availability, the hospital rarely operates more than two simultaneous rooms for the same specialty. Therefore, the following capacity constraints were adopted throughout the experiments:
\[
N_{mt}^{\rm AM}
=
N_{mt}^{\rm PM}
=
2,
\qquad
\forall m\in M,\ t\in T.
\]

\subsection{Surgical duration estimation}

To estimate surgery durations, we analyzed a database containing electronic records of 12,054 surgical interventions collected during 26 consecutive months of hospital activity.

Each registry contains:
\begin{itemize}
    \item patient information;
    \item surgical specialty;
    \item intervention date;
    \item surgical-team information;
    \item and timestamps associated with four stages of the intervention:
    \begin{enumerate}
        \item pre-anesthesia,
        \item anesthesia,
        \item surgical act,
        \item recovery.
    \end{enumerate}
\end{itemize}

\begin{figure*}[h]
\centering
\includegraphics[width = 16cm]{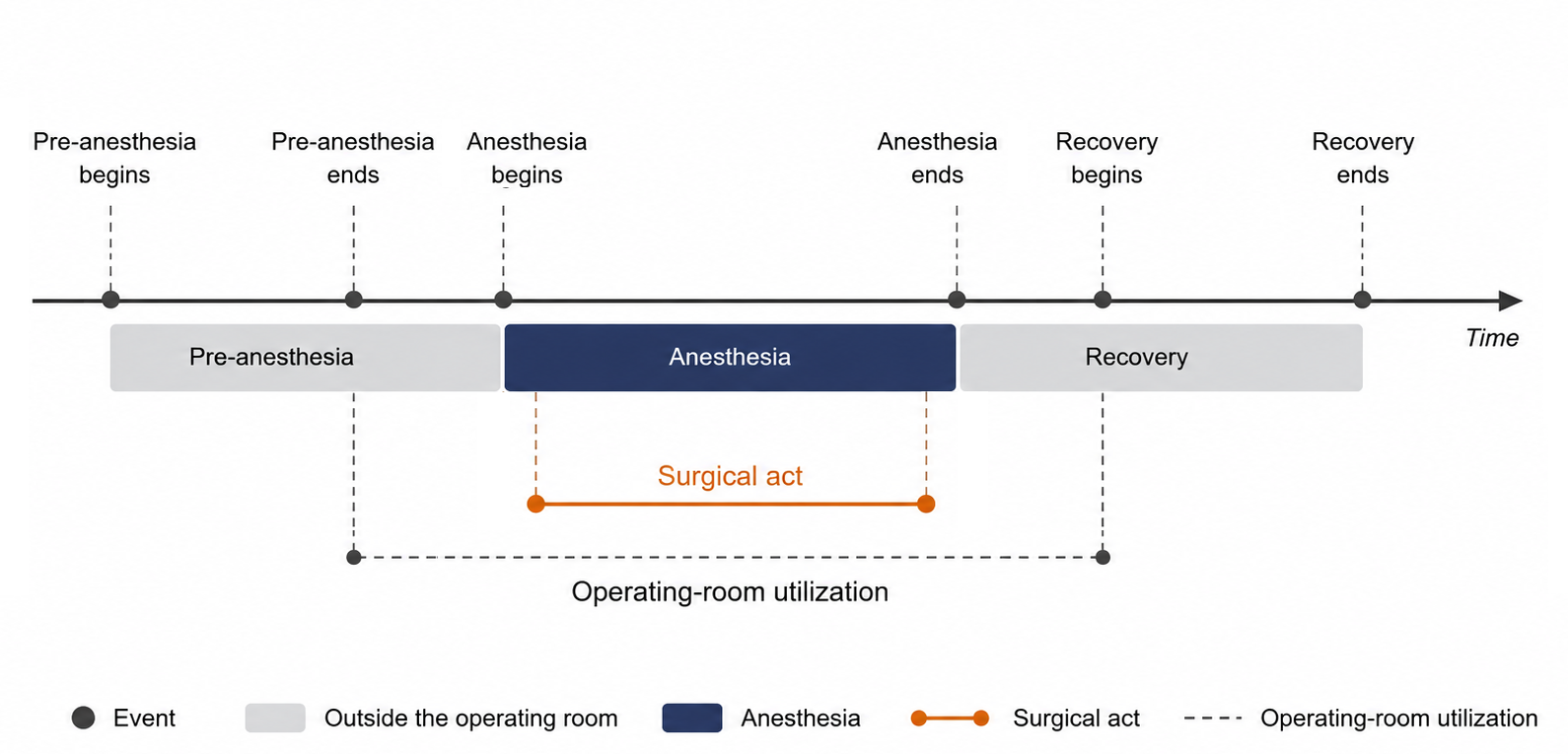}
\caption{Timeline of the surgical process and operating-room utilization during an elective surgical intervention.}
\label{fig:surgical_process_timeline}
\end{figure*}

The operating room is occupied during the anesthesia stage, which includes the surgical act, together with preparation and post-procedure activities, as illustrated in Figure~\ref{fig:surgical_process_timeline}.

Each intervention may contain multiple surgical procedures. Consequently, the exact operating-room occupation time associated with a single procedure is generally unavailable. However, the database identifies a principal surgery for each intervention. Therefore, the duration associated with a surgical procedure was estimated as the average operating-room utilization time among interventions where that procedure appeared as the principal operation.

Although this approximation does not explicitly model patient-level variability in surgery duration, both compared methodologies (\block\ and \milp) use the same duration estimates. Thus, the experimental comparison isolates the effect of the allocation methodology itself. Improving surgery-duration prediction constitutes a relevant direction for future work and could naturally be integrated into the proposed framework.

\subsection{Initial waiting-list characterization}

Using the estimated surgery durations, we computed the demanded surgical time associated with the waiting list of each specialty.

\begin{figure*}[h]
\centering
\includegraphics[width = 2\columnwidth]{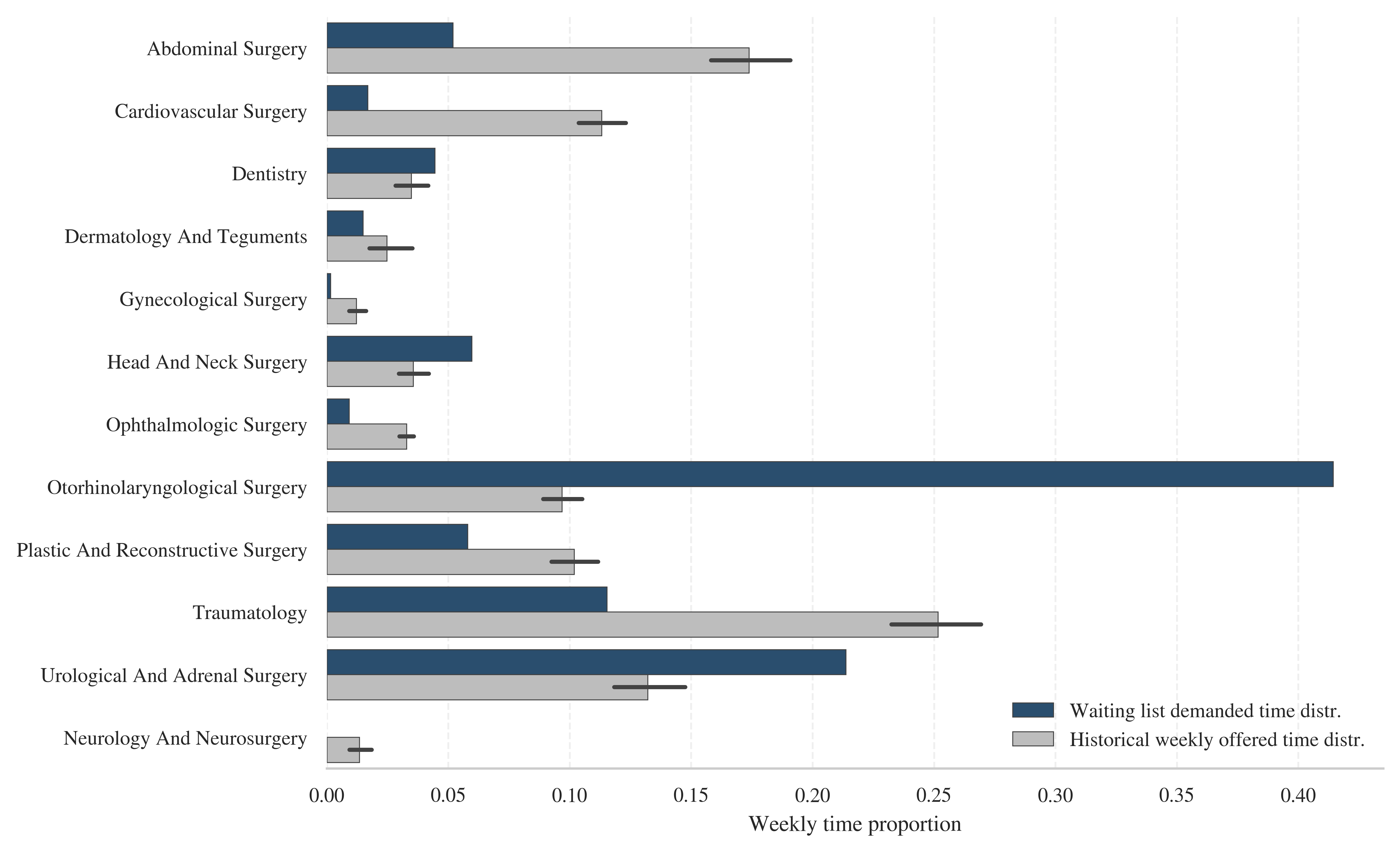}
\caption{Comparison between the historical proportion of weekly operating-room time offered by the hospital and the estimated proportion of demanded surgical time in the waiting list, grouped by medical specialty.}
\label{fig:historicalTimeDistribution}
\end{figure*}

Figure~\ref{fig:historicalTimeDistribution} compares:
\begin{itemize}
    \item the historical proportion of weekly offered surgical time by specialty,
    \item and the estimated demanded surgical time derived from the waiting list.
\end{itemize}

The comparison reveals an important imbalance between supply and demand across specialties. 
In particular, Otorhinolaryngological Surgery and Urology exhibit substantially larger waiting-list demand than the historically allocated operating-room capacity, whereas other specialties appear relatively over-served.

Table~\ref{tab:statWaitingList} summarizes the initial state of the surgical waiting list used in the simulation experiments, including:
\begin{itemize}
    \item number of patients by specialty;
    \item mean waiting time;
    \item median waiting time;
    \item maximum waiting time;
    \item and standard deviation.
\end{itemize}

\begin{table*}
\centering
\begin{tabular}{lrrrrr}
\toprule
\multirow{2}{*}{Medical specialty}& List size &\multicolumn{4}{c}{Waiting time stats (days in list)}\\
 &	(\# patients) &	mean &	median &	max &	std\\
\midrule
Abdominal surgery &	28&	187&	130&	738&	197.4\\
Cardiovascular surgery&	7&	186&	36&	501&	228.9\\
Dentistry&	33&	144&	96&	537&	139.4\\
Dermatology and teguments&	20&	329&	298&	799&	216.1\\
Gynecological surgery&	1&	33&	33&	33 & \\	
Head and neck surgery&	22&	285&	274&	579&	182.8\\
Ophthalmologic surgery&	5&	116&	80&	269&	89.5\\
Otorhinolaryngological surgery&	636&	243&	221&	921&	176.5\\
Plastic and reconstructive surgery&	35&	225&	171&	838&	218.1\\
Traumatology&	53&	343&	403&	905&	286.0\\
Urological and adrenal surgery&	190&	143&	91&	780&	170.9\\ \midrule
Full waiting list&	1030&	226&	171&	921&	191.3\\
\bottomrule
\end{tabular}
\caption{Summary statistics of the initial waiting list used in the simulation experiments, including waiting times (in days) by medical specialty.}
\label{tab:statWaitingList}
\end{table*}

These statistics provide the initial conditions for the simulation experiments and illustrate the operational pressure experienced by several specialties.

\subsection{Performance metrics}

The proposed framework (\block) was compared against the benchmark \milp\ formulation using the following performance indicators:

\begin{itemize}
    \item \textbf{Patients operated (PO):} average number of surgeries successfully performed;
    
    \item \textbf{Mean waiting time:} average waiting time of patients remaining in the waiting list;
    
    \item \textbf{Median waiting time:} median waiting time of patients remaining in the waiting list;
    
    \item \textbf{Utilization rate (UR):} proportion of operating-room occupation time relative to the available operating-room capacity;
    
    \item \textbf{Offered/demanded time balance:} comparison between assigned operating-room time and estimated demanded time by specialty.
\end{itemize}

Additionally, the long-term experiments incorporate fairness analyses using Lorenz curves and Gini coefficients to evaluate the equity of operating-room time distribution among specialties.

\subsection{Simulation assumptions}

Two operational phenomena observed in hospital practice were incorporated into the simulations.

\subsubsection{Patient confirmation process}

Before scheduling a surgery, hospital personnel must contact the patient and obtain confirmation of the proposed intervention date. Historical records indicate that approximately one out of four contacted patients confirms the appointment.

To model this phenomenon, each patient $p\in P$ was assigned an independent Bernoulli random variable:
\[
B_p \sim \text{Bernoulli}(0.25),
\]
where:
\[
B_p = 1
\]
indicates that patient $p$ successfully confirms the proposed surgery date.

This confirmation mechanism was incorporated into the patient-allocation procedure during the simulations.

\subsubsection{Surgery suspension process}

Even after confirmation, surgeries may be suspended due to medical or administrative reasons, including illness, incomplete examinations, or operational contingencies.

Historical data indicate a surgery suspension rate close to $7\%$. Therefore, after patient assignment, scheduled surgeries were independently canceled according to a Bernoulli process with parameter:
\[
p = 0.07.
\]

These suspended surgeries were accounted for when computing the final performance indicators.

\subsection{Simulation scenarios}

Two experimental scenarios were considered.

\subsubsection{Single-period scheduling}

The first scenario evaluates a one-week scheduling process under high-demand conditions. The  experiment considers:
\begin{itemize}
    \item five operating rooms;
    \item one planning horizon of five working days;
    \item and waiting-time prioritization (FIFO criterion).
\end{itemize}

The objective of this scenario is to compare the short-term operational behavior of \block\ and \milp\ in terms of:
\begin{itemize}
    \item operating-room utilization;
    \item waiting-time reduction;
    \item computational time;
    \item and specialty-demand balancing.
\end{itemize}

\subsubsection{Multi-period scheduling}

The second scenario evaluates the long-term behavior of the proposed framework over a 20-week horizon under a weekly replanning regime.

The simulations assume:
\begin{itemize}
    \item availability of three operating rooms;
    \item fixed operational conditions throughout the experiment;
    \item and weekly arrivals of new patients.
\end{itemize}

Historical analyses of patient admissions showed that the weekly number of new patients entering the waiting list follows an approximately Normal distribution, as illustrated in Figure~\ref{fig:admissionPatients}. Therefore, weekly arrivals were simulated as:
\[
A_t \sim \mathcal{N}(37,11.8^2).
\]

The surgical specialties associated with newly arriving patients were generated according to historical proportions estimated from 222 different surgical procedures, summarized in Table~\ref{tab:proporciones_operaciones}.

The purpose of this second scenario is to evaluate the ability of the proposed framework to sustain long-term reductions in waiting times while preserving balanced specialty allocation and high operating-room utilization.

\section{Results}\label{sec:results}

The proposed two-stage framework (\block) is compared against the benchmark patient-level optimization approach (\milp) in terms of:
\begin{itemize}
    \item waiting-time reduction;
    \item operating-room utilization;
    \item specialty-demand balancing;
    \item computational performance;
    \item and long-term operational behavior.
\end{itemize}

\subsection{Single-period scheduling results}

We first evaluate the behavior of the models under a single weekly planning horizon.

Table~\ref{tab:MILPvsBLOCK_general} summarizes the main performance indicators obtained by both methodologies. The benchmark \milp\ model achieved a slightly larger number of operated patients and higher operating-room utilization than \block. In particular, the utilization ratio obtained by \milp\ reached $98.1\%$, whereas \block\ achieved a utilization ratio of $90.9\%$.

\begin{table*}[ht]
    \centering
    \begin{tabular}{|c|c|c|c}
    \hline
    Key Performance Indicator                            & {\milp}  & {\block}\\ \hline
    Number of operations              & 144   & 130   \\
    Final average waiting time (days) & 224.1 & 225.6 \\
    Final median waiting time (days)    & 200   & 204   \\
    Utilization ratio (\%)            & 98.1  & 90.9  \\ 
    Computation time (secs)              & 1,800      & 1.6      \\
    \hline
    \end{tabular}
    \caption{Comparison of the main performance indicators obtained by the \milp\ and \block\ models under the single-period scheduling scenario.}
    \label{tab:MILPvsBLOCK_general}
\end{table*}


The computational performance of the two methodologies also differs substantially. While the \milp\ formulation required approximately $1{,}800$ seconds to obtain a solution, the proposed \block\ framework required only $1.6$ seconds. This difference is particularly relevant in practical hospital settings, where weekly replanning must be performed repeatedly under changing waiting-list conditions.

Moreover, the proposed framework preserves the proportional relationship between demanded and offered specialty time more accurately than the benchmark formulation. Table~\ref{tab:outputBLOCK} shows that the operating-room time assigned by \block\ closely follows the demanded surgical time for most specialties. In contrast, \milp\ tends to under-allocate several specialties despite its higher overall operating-room utilization.

\begin{table*}
\centering
\resizebox{16cm}{!} {
\begin{tabular}{lrrrrrr}
\toprule
\multirow{2}{*}{Specialty}& \multirow{2}{*}{Demanded hours} &\multicolumn{2}{c}{Offered hours by models} &\multicolumn{3}{c}{Number of blocks}\\
& &	{\milp} & {\block} &	AM &	PM &	Friday PM\\
\midrule
Urology and Nephrology         &             44.3 &          43.0 &            43 &                 6 &                         3 &                        2 \\
Dentistry                      &             11.9 &           9.8 &            13 &                 2 &                         1 &                        0 \\
Gastroenterology               &              5.7 &           4.1 &             5 &                 1 &                         0 &                        0 \\
Traumatology                   &             38.9 &          37.4 &            39 &                 5 &                         4 &                        1 \\
Otorhinolaryngological Surgery &             73.3 &          80.8 &            73 &                 9 &                         8 &                        2 \\
Head and Neck Surgery          &                2.0 &           0.0 &             3 &                 0 &                         1 &                        0 \\
Dermatology and Teguments      &              2.4 &           2.4 &             3 &                 0 &                         1 &                        0 \\
Plastic Surgery                &             14.3 &          12.0 &            13 &                 2 &                         1 &                        0 \\
Cardiovascular Surgery         &              2.3 &           0.0 &             3 &                 0 &                         1 &                        0 \\
Neurology and Neurosurgery     &                0.0 &           0.0 &             0 &                 0 &                         0 &                        0 \\
\bottomrule
\end{tabular}}
\caption{Comparison between the operating-room time offered by the \milp\ and \block\ models and the estimated demanded surgical time across medical specialties.}
\label{tab:outputBLOCK}
\end{table*}

In particular, specialties such as Cardiovascular Surgery and Head and Neck Surgery receive no allocated time under the \milp\ formulation, whereas \block\ guarantees at least a minimum operational allocation consistent with the balancing constraints imposed by the tactical block-distribution stage.

Table~\ref{tab:ORuse} presents the operating-room utilization generated by \block\ during the weekly planning horizon. The proposed framework maintains high operating-room occupation levels across all operating rooms, achieving an overall weekly utilization ratio of $91.8\%$, while preserving operational flexibility and specialty balancing.

\begin{table*}[ht]
\centering
\setlength{\tabcolsep}{5pt}
\renewcommand{\arraystretch}{1.15}

\begin{tabular}{lcccccc}
\toprule
& \multicolumn{5}{c}{Operating room} & \multirow{2}{*}{Daily UR (\%)} \\
\cmidrule(lr){2-6}
Day & Room 1 & Room 2 & Room 3 & Room 4 & Room 5 &  \\

\midrule
Monday    & 89.5 & 89.7 & 96.6 & 97.2 & 88.7 & 92.3 \\
Tuesday   & 96.4 & 94.3 & 98.6 & 88.1 & 69.6 & 87.4 \\
Wednesday & 82.6 & 89.5 & 92.5 & 95.4 & 89.9 & 90.0 \\
Thursday  & 87.3 & 88.3 & 90.9 & 95.4 & 90.9 & 90.5 \\
Friday    & 76.8 & 83.0 & 71.7 & 93.8 & 94.9 & 84.0 \\

\midrule
Weekly UR (\%) & 89.5 & 92.0 & 93.4 & 97.0 & 87.3 & \textbf{91.8} \\

\bottomrule
\end{tabular}
\caption{Operating-room utilization rates generated by the \block\ framework during the single-week planning horizon, grouped by operating room and day.}
    \label{tab:ORuse}
\end{table*}

These results suggest that the tactical decomposition introduced by \block\ produces a more balanced allocation of surgical capacity without severely compromising utilization levels.

\subsection{Multi-period scheduling results}

We next evaluate the long-term behavior of the scheduling methodologies over a rolling 20-week planning horizon.

Table~\ref{tab:MILPvsBLOCK_weeks} summarizes the aggregated results obtained throughout the simulation period. In contrast to the single-period scenario, the proposed \block\ framework outperforms \milp\ in terms of waiting-list reduction and number of operated patients.

\begin{table*}[ht]
    \centering
    \begin{tabular}{|c|c|c|c}
    \hline
    Key Performance Indicator                            & {\milp}  & {\block}\\ \hline
    Number of operations              & 1,207   & 1,337   \\
    Final average waiting time (days) & 78.1 & 60.3 \\
    Final median waiting time (days)    & 70   & 42   \\
    Average utilization ratio (\%)            & 92.1  & 84  \\ 
    Computation time (secs)              & 21,600      & 32      \\
    \hline
    \end{tabular}
    \caption{Comparison of the long-term performance of {\milp} and {\block} over the 20-week rolling-horizon simulation.}
    \label{tab:MILPvsBLOCK_weeks}
\end{table*}

More precisely, \block\ performs a total of $1{,}337$ surgeries during the simulation period, compared to $1{,}207$ surgeries under \milp. Furthermore, the final mean waiting time is reduced from $78.1$ days under \milp\ to $60.3$ days under \block, corresponding to an approximate reduction of $23\%$.

The reduction in median waiting time is even more pronounced. The median waiting time decreases from $70$ days under \milp\ to $42$ days under \block, corresponding to an approximate reduction of $40\%$. These results indicate that the proposed framework not only improves average system performance, but also reduces the waiting-time burden experienced by a large proportion of patients.


The two methodologies exhibit different operational characteristics with respect to operating-room utilization. The \milp\ model achieves an average utilization ratio of $92.1\%$, whereas the proposed \block\ framework maintains a utilization ratio close to $84\%$. As shown previously, this reduction in utilization is accompanied by substantial improvements in waiting-list dynamics, specialty balancing, and computational efficiency.

The computational advantage of the proposed framework also becomes more evident in the long-term simulations. The total computational time required by \milp\ exceeded $21{,}600$ seconds, whereas \block\ required only $32$ seconds over the complete 20-week simulation horizon.

\begin{figure*}[htbp]
    \centering

    \begin{subfigure}{\linewidth}
        \centering
        \includegraphics[width=0.8\linewidth]{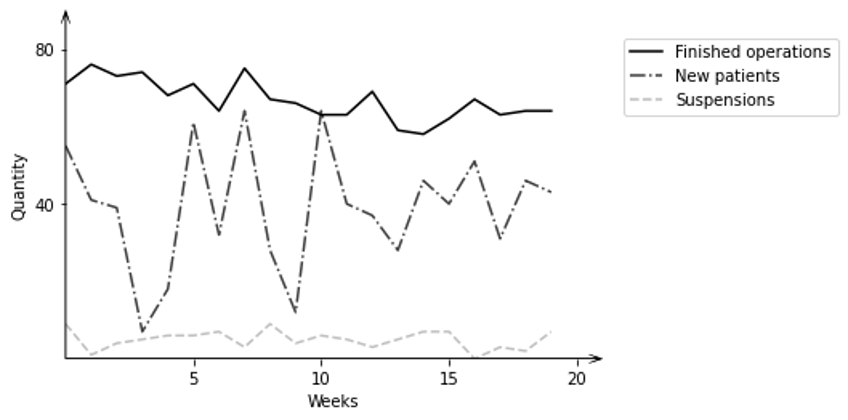}
        \caption{{\block}}
        \label{fig:operations_block}
    \end{subfigure}
    \hfill
    \begin{subfigure}{\linewidth}
        \centering
        \includegraphics[width=0.8\linewidth]{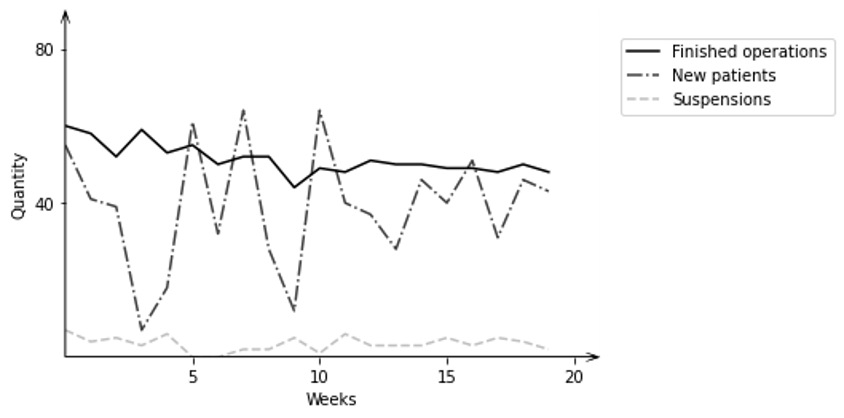}
        \caption{{\milp}}
        \label{fig:operations_milp}
    \end{subfigure}

    \caption{Weekly evolution of completed surgeries, incoming patients, and suspended surgeries over the 20-week simulation period for (a) {\block} and (b) {\milp}.}
    \label{fig:patients_weeks}
\end{figure*}

Figure~\ref{fig:patients_weeks} illustrates the weekly evolution of:
\begin{itemize}
    \item completed surgeries;
    \item incoming patients;
    \item and surgery suspensions.
\end{itemize}

The proposed framework maintains a relatively stable number of completed surgeries throughout the simulation period despite fluctuations in patient arrivals and surgery suspensions.

\begin{figure*}[htbp]
    \centering

    \begin{subfigure}{\linewidth}
        \centering
        \includegraphics[width=0.6\linewidth]{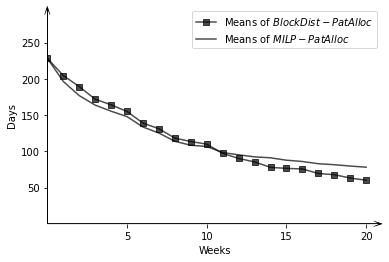}
        \caption{}
        \label{fig:mean_b_m}
    \end{subfigure}
    \hfill
    \begin{subfigure}{\linewidth}
        \centering
        \includegraphics[width=0.6\linewidth]{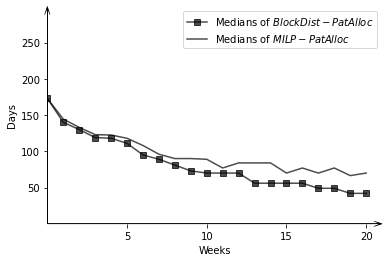}
        \caption{}
        \label{fig:median_b_m}
    \end{subfigure}

    \caption{Weekly evolution of (a) mean waiting time and (b) median waiting time over the 20-week simulation period for {\block} and {\milp}.}
    \label{fig:median_mean}
\end{figure*}

Figure~\ref{fig:median_mean} shows the weekly evolution of the mean and median waiting times over the 20-week simulation period for both methodologies. Although both approaches initially reduce waiting times, the reductions achieved by \block\ become progressively larger over time. This behavior is particularly evident for the median waiting time, where the proposed framework produces a substantially faster decrease during the second half of the simulation horizon.

\subsection{Utilization and fairness analysis}

\begin{figure*}[h]
    \centering
    \begin{subfigure}{0.7\columnwidth}
        \centering
        \includegraphics[width=1.2\linewidth]{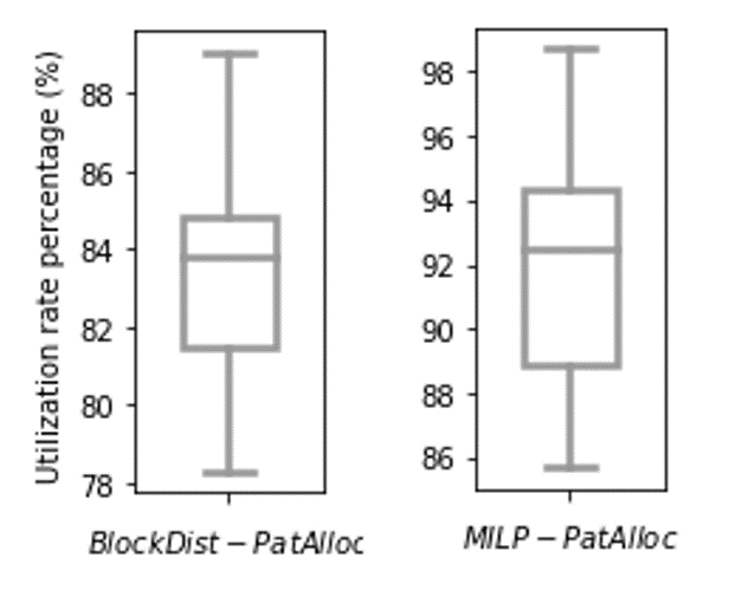}
        \caption{}
    \end{subfigure}
    \hfill
    \begin{subfigure}{1.2\columnwidth}
        \centering
        \includegraphics[width=\linewidth]{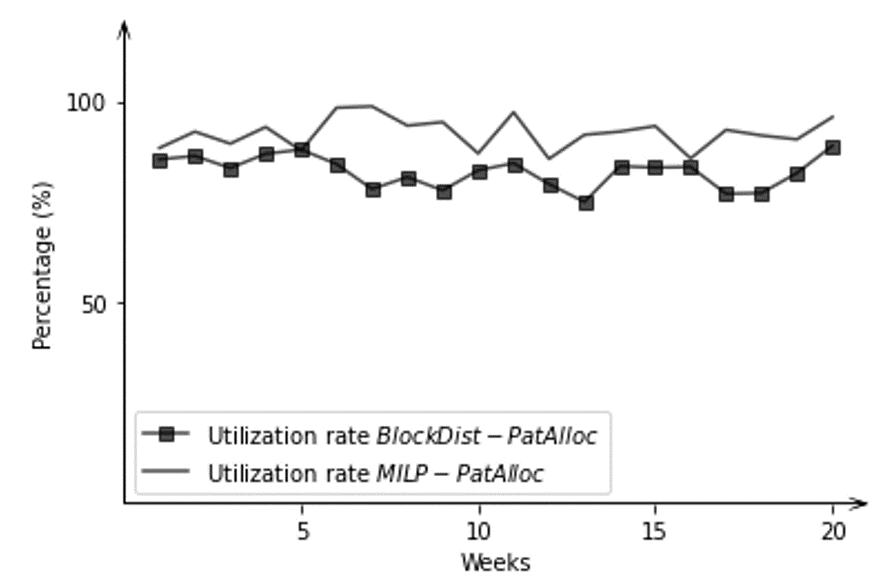}
        \caption{}
    \end{subfigure}
    \caption{Utilization rate over the 20-week simulation period. The left panel (a) shows the distribution of weekly utilization rates for  {\block} and  {\milp}, while the right panel (b) shows the week-by-week evolution of the utilization rate for both methods.}
    \label{fig:ur_weeks}
\end{figure*}

Figure~\ref{fig:ur_weeks} shows the weekly evolution of the operating-room utilization rates during the 20-week simulations. As expected, \milp\ systematically achieves higher utilization levels, often close to full operating-room occupation. Nevertheless, the proposed \block\ framework maintains relatively stable utilization values while simultaneously producing better waiting-time reductions.

This observation suggests that maximizing operating-room occupation alone does not necessarily translate into improved long-term waiting-list performance. Instead, the balanced allocation of operating-room capacity across specialties appears to play a central role in reducing waiting-list pressure.

To further evaluate the equity of the proposed allocations, Figure~\ref{fig:lorenzGini} presents the Lorenz curves associated with the offered-demanded specialty distributions for both methodologies.
\begingroup \colorlet{mygray}{black!65}
\begin{figure*}[htbp]
    \centering

    \begin{subfigure}{\columnwidth}
        \centering
        \includegraphics[width=\linewidth]{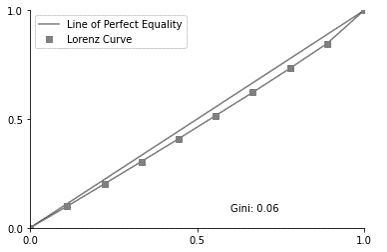}
        \caption{{\block}}
        \label{fig:lorenzGini_block}
    \end{subfigure}
    \hfill
    \begin{subfigure}{\columnwidth}
        \centering
        \includegraphics[width=\linewidth]{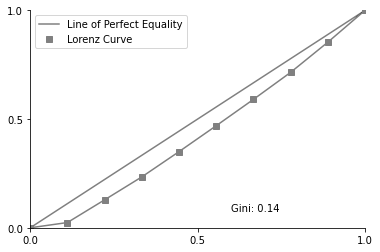}
        \caption{{\milp}}
        \label{fig:lorenzGini_milp}
    \end{subfigure}

    \caption{Lorenz curves associated with the distribution of offered operating-room time relative to demanded surgical time across specialties for (a) {\block} and (b)   {\milp}. The solid diagonal line represents perfect equality, whereas the square-marked curve represents the observed distribution of operating-room allocation. Lower deviations from the equality line indicate a more balanced distribution among specialties.}
    \label{fig:lorenzGini}
\end{figure*}
\endgroup

The proposed \block\ framework achieves a substantially lower Gini coefficient ($0.06$) compared to \milp\ ($0.14$), indicating a significantly more equitable distribution of operating-room time among specialties relative to their demanded surgical time.

This result highlights one of the main strengths of the proposed two-stage framework: the explicit incorporation of specialty-demand balancing during the tactical allocation stage leads to more homogeneous access to operating-room capacity across specialties.

Finally, Table~\ref{tab:perc_hours_week} confirms that the offered-to-demanded time ratios generated by \block\ remain close to one for nearly all specialties throughout the simulation horizon. In contrast, the benchmark \milp\ formulation systematically under-serves several specialties despite achieving higher operating-room utilization.

\begin{table*}[htbp]
\centering
\resizebox{2\columnwidth}{!}{
\begin{tabular}{lrrr|rrr}
\toprule
& \multicolumn{3}{c|}{\milp} & \multicolumn{3}{c}{\block} \\
\cmidrule(lr){2-4} \cmidrule(lr){5-7}
{} & Req Hours & Off Hours & Off/Req Time &
Req Hours & Off Hours & Off/Req Time \\
\midrule

Plastic surgery & 156.94 & 94.92 & 0.60 & 157.70 & 146.00 & 0.93 \\
Urology and Nephrology & 588.51 & 365.43 & 0.62 & 571.03 & 568.00 & 0.99 \\
Traumatology & 250.50 & 149.08 & 0.60 & 230.29 & 219.00 & 0.95 \\
Otorhinolaryngology surgery & 823.93 & 602.82 & 0.73 & 831.20 & 833.00 & 1.00 \\
Cardiovascular surgery & 48.59 & 5.87 & 0.12 & 49.20 & 52.00 & 1.06 \\
Dentistry & 167.29 & 87.22 & 0.52 & 149.03 & 153.00 & 1.03 \\
Gastroenterology & 130.23 & 75.55 & 0.58 & 134.91 & 134.00 & 0.99 \\
Head and neck surgery & 119.17 & 62.77 & 0.53 & 143.46 & 140.00 & 0.98 \\
Dermatology & 54.83 & 38.00 & 0.69 & 70.90 & 70.00 & 0.99 \\

\bottomrule
\end{tabular}
}
\caption{Comparison between the total demanded surgical time and the total operating-room time assigned by each methodology during the 20-week simulation period, grouped by medical specialty.}
\label{tab:perc_hours_week}
\end{table*}

Overall, the results indicate that the proposed \block\ framework provides:
\begin{itemize}
    \item substantial reductions in mean and median waiting times;
    \item improved specialty-level fairness;
    \item stable operating-room utilization;
    \item and significantly lower computational requirements.
\end{itemize}

\section{Discussion}\label{sec:discussion}

The results obtained in both experimental scenarios indicate that the proposed two-stage framework (\block) constitutes an effective and computationally efficient approach for operating-room allocation under persistent waiting-list pressure in public hospitals.

A first important observation is that maximizing operating-room utilization alone does not necessarily lead to improved waiting-list performance. Although the benchmark \milp\ formulation systematically achieved higher utilization rates, the proposed \block\ framework consistently produced lower mean and median waiting times during the long-term simulations. 
This behavior suggests that the distribution of operating-room capacity across specialties plays a critical role in the evolution of waiting lists.

The proposed framework explicitly incorporates this balancing mechanism through the tactical block-allocation stage. By assigning operating-room capacity proportionally to the estimated surgical demand of each specialty, the framework avoids the systematic under-service of specialties with high accumulated backlog. This effect becomes particularly visible in Table~\ref{tab:perc_hours_week}, where the offered-to-demanded time ratios generated by \block\ remain close to one for most specialties throughout the simulation horizon.

In contrast, the benchmark \milp\ formulation tends to concentrate operating-room usage on specialties that maximize short-term operational efficiency, producing higher utilization levels but generating larger inequalities in specialty access. This phenomenon is also reflected in the Lorenz curves and Gini coefficients shown in Figure~\ref{fig:lorenzGini}. 
The substantially lower Gini coefficient obtained by \block\ indicates that the proposed framework produces a more equitable distribution of operating-room capacity relative to specialty demand.

From an operational perspective, this result is particularly relevant in public healthcare systems, where waiting-list management must balance efficiency and fairness simultaneously. 
Focusing exclusively on operating-room occupation may improve short-term throughput while progressively increasing inequities between specialties. The proposed framework partially mitigates this phenomenon by explicitly incorporating specialty-demand balancing into the tactical scheduling stage.

Another important result is the substantial computational advantage achieved by the proposed framework. The decomposition into tactical and operational stages dramatically reduces computational complexity compared to the direct patient-level MILP formulation. In the 20-week simulations, the benchmark formulation required approximately $21{,}600$ seconds of computational time, whereas \block\ required only $32$ seconds.

This computational reduction is highly relevant in practice. Hospital waiting lists evolve continuously due to:
\begin{itemize}
    \item arrival of new patients;
    \item emergency prioritizations;
    \item surgery suspensions;
    \item and unsuccessful patient confirmations.
\end{itemize}

Consequently, operating-room schedules must be recomputed frequently, typically every week. 
Under these conditions, computationally expensive patient-level optimization models become difficult to deploy operationally, particularly in hospitals with limited technological infrastructure.

The decomposition adopted in \block\ also introduces an important degree of operational flexibility. Since patient assignment is separated from specialty block allocation, modifications in the patient schedule can be performed without recomputing the complete tactical allocation. This property becomes particularly valuable in the presence of cancellations or failed confirmations, which occur frequently in real hospital environments.

Indeed, the simulations explicitly incorporated two operational phenomena observed in practice:
\begin{itemize}
    \item low patient-confirmation rates;
    \item and surgery suspensions.
\end{itemize}

The proposed framework remained stable under these perturbations, maintaining high utilization levels and sustained reductions in waiting times throughout the simulation horizon.

The experiments also highlight the importance of medium-term planning. While the differences between both methodologies are relatively moderate during a single scheduling period, the long-term simulations reveal cumulative effects associated with the balancing mechanism introduced by \block. Over time, these effects produce substantial reductions in waiting-list pressure, particularly in the median waiting time.

This observation suggests that tactical balancing decisions may have a stronger impact on long-term waiting-list dynamics than purely local patient-level optimization criteria.

The present work nevertheless has several limitations. First, surgery durations were estimated using historical average intervention times associated with the principal surgical procedure. Although this assumption allows a fair comparison between methodologies, it does not capture patient-level variability in surgical duration.

Second, the experiments consider a simplified representation of patient prioritization based primarily on waiting time (FIFO criterion). In practice, hospitals frequently incorporate additional clinical considerations, including urgency and diagnosis severity. However, the proposed framework naturally allows the incorporation of arbitrary hospital-defined priority scores.

Finally, the current framework does not explicitly incorporate emergency surgeries, stochastic operating-room disruptions, or dynamic rescheduling during the execution of the surgical day. 
These aspects constitute relevant directions for future research.

Despite these limitations, the results indicate that the proposed two-stage framework provides an effective compromise between:
\begin{itemize}
    \item computational tractability;
    \item operational flexibility;
    \item specialty-level fairness;
    \item and waiting-list reduction.
\end{itemize}

Therefore, the proposed methodology appears particularly well suited for public hospital environments, where operating-room allocation decisions must simultaneously address resource limitations, dynamic waiting lists, and heterogeneous specialty demands.

\section{Conclusion}\label{sec:conclusion}

In this work, we proposed a two-stage operating-room allocation framework designed to support surgical scheduling decisions in public hospital environments under persistent waiting-list pressure. The proposed methodology combines:
\begin{itemize}
    \item a tactical operating-room block allocation stage (\blockDist),
    \item and an operational patient scheduling procedure (\patAlloc).
\end{itemize}

The framework explicitly incorporates specialty-demand balancing into the allocation process while preserving operational flexibility and computational tractability.

The experimental results obtained using historical data from Dr. Luis Calvo Mackenna Hospital indicate that the proposed framework substantially improves long-term waiting-list performance compared to a direct patient-level optimization benchmark. In particular, the proposed methodology achieved:
\begin{itemize}
    \item significant reductions in mean and median waiting times;
    \item a larger number of completed surgeries over the simulation horizon;
    \item improved balance between offered and demanded specialty capacity;
    \item and substantially lower computational times.
\end{itemize}


The benchmark patient-level MILP formulation achieved higher operating-room utilization, whereas the proposed framework consistently produced better waiting-list dynamics and a more equitable distribution of operating-room time across specialties. These findings highlight that operating-room utilization, while an important performance indicator, should be considered jointly with waiting-list reduction and equitable resource allocation when evaluating scheduling methodologies for public healthcare systems.

One of the main strengths of the proposed methodology lies in the hierarchical decomposition between tactical and operational decisions. This decomposition considerably reduces computational complexity while simultaneously introducing important operational flexibility. In practice, this allows hospitals to recompute schedules frequently under changing waiting-list conditions without requiring computationally expensive patient-level optimization procedures.

Another important contribution of the framework is its adaptability to real hospital environments. The proposed methodology naturally supports:
\begin{itemize}
    \item rolling-horizon weekly scheduling;
    \item specialty-specific operational constraints;
    \item patient prioritization policies;
    \item surgery suspensions;
    \item and failed patient confirmations.
\end{itemize}

These characteristics make the framework particularly suitable for public hospital systems operating under dynamic demand conditions and limited resource availability.

The results suggest that the proposed two-stage framework constitutes a practical and scalable alternative for operating-room allocation in public healthcare systems, providing an effective balance between operational efficiency, fairness across specialties, and long-term waiting-list reduction.

\section*{Declarations}

\bmhead{\bf Funding}
This work was supported by the Agencia Nacional de Investigación y Desarrollo (ANID) through ANID BASAL FB210005, which supports the Centro de Modelamiento Matemático (CMM) as a Center of Excellence.

\bmhead{\bf Acknowledgements}
The authors thank Hospital Dr. Luis Calvo Mackenna (HLCM) for providing the anonymized operational data used in this study and for sharing valuable institutional knowledge regarding elective surgery scheduling and waiting-list management.

\bmhead{\bf Conflict of interest}
The authors declare that they have no competing interests.

\bmhead{\bf Ethics approval and consent to participate}
This study used anonymized operational records provided by Hospital Dr. Luis Calvo Mackenna. The data did not contain patient-identifiable information. Therefore, the study did not involve direct interaction with patients and did not require informed consent.

\bmhead{\bf Data availability}
The data analyzed in this study were provided by Hospital Dr. Luis Calvo Mackenna under institutional conditions and are not publicly available. Aggregated or processed data may be made available from the corresponding author upon reasonable request and subject to institutional approval.

\bmhead{\bf Code availability}
The code used for the implementation of the optimization models and numerical experiments may  be made available from the corresponding author upon reasonable request.

\bmhead{\bf Author contribution}
All authors contributed to the conception and design of the study. Data analysis, model development, computational implementation, and interpretation of results were performed collaboratively by the authors. All authors contributed to the writing, revision, and approval of the final manuscript.

\bmhead{\bf Declaration of Generative AI and AI-Assisted Technologies in the Writing Process}
During the preparation of this work, the authors used AI-assisted writing tools to improve the grammar, clarity, and readability of the manuscript. All suggestions generated by these tools were carefully reviewed, edited, and verified by the authors, who take full responsibility for the final content of the article.

\begin{appendices}

\section{Patient-level MILP benchmark formulation (\milp)}\label{MILP}

In this work, we compare the proposed two-stage framework (\block) with a patient-level mixed-integer linear programming formulation, which we denote by \milp. This benchmark model is inspired by the formulation proposed in \cite{wolff2012}, with two main modifications.

First, instead of maximizing operating-room occupation while penalizing overtime, we prioritize patients according to a predefined priority criterion. Second, we incorporate an explicit balancing mechanism between the offered operating-room time and the demanded surgical time across specialties.

The formulation directly assigns patients to operating rooms and planning days, without introducing a tactical specialty-block allocation stage.

Using the notation introduced in Table~\ref{tabla_constantes}, we define the additional parameters and variables required by the benchmark formulation.

\subsection*{Additional parameter}

\begin{itemize}
    \item $c_{rt}$: maximum operating-room time capacity available in room $r \in R$ during day $t \in T$.
\end{itemize}

\subsection*{Decision variables}

\begin{itemize}
    \item $x_{prt}$: binary variable equal to $1$ if patient $p \in P$ is assigned to operating room $r \in R$ on day $t \in T$, and $0$ otherwise.
    
    \item $g_m$: deviation between the operating-room time assigned to specialty $m \in M$ and its ideal proportional demanded time.
    
    \item $\tilde g_m$: auxiliary nonnegative variable used to linearize the absolute value of $g_m$.
\end{itemize}

\subsection*{Auxiliary quantities}

The total operating-room time assigned to specialty $m$ is defined as:
\begin{equation*}
Y_m =
\sum_{p \in P_m}
\sum_{r \in R}
\sum_{t \in T}
x_{prt}(d_p + u),
\quad m \in M.
\end{equation*}

The total assigned operating-room time is:
\begin{equation*}
Y =
\sum_{p \in P}
\sum_{r \in R}
\sum_{t \in T}
x_{prt}(d_p + u)
=
\sum_{m \in M} Y_m.
\end{equation*}

The demanded surgical times by specialty ($D_m$) and globally ($D$) are defined identically to those introduced in the \blockDist\ formulation.

\subsection*{Constraints}

The benchmark formulation considers the following constraints.

\begin{enumerate}

\item \textbf{Operating-room capacity constraints}

The total assigned surgical and cleaning time in each operating room and day cannot exceed the available room capacity:
\begin{equation}
\sum_{p \in P}
x_{prt}(d_p + u)
\leq
c_{rt},
\quad
\forall r \in R,\ t \in T.
\end{equation}

\item \textbf{Single-assignment constraints}

Each patient can be assigned at most once within the planning horizon:
\begin{equation}
\sum_{r \in R}
\sum_{t \in T}
x_{prt}
\leq 1,
\quad
\forall p \in P.
\end{equation}

\item \textbf{Demand-balance constraints}

The variable $g_m$ measures the deviation between the assigned operating-room time and the ideal proportional allocation associated with specialty demand:
\begin{equation}
Y_m
=
g_m
+
\frac{D_m}{D}Y,
\quad
\forall m \in M.
\end{equation}

\item \textbf{Absolute-value linearization}

The auxiliary variables $\tilde g_m$ are used to linearize the absolute values of the deviations:
\begin{equation}
g_m \leq \tilde g_m,
\qquad
-g_m \leq \tilde g_m,
\quad
\forall m \in M.
\end{equation}

\end{enumerate}

\subsection*{Objective function}

The benchmark formulation maximizes patient prioritization while penalizing imbalances between offered and demanded specialty capacity:
\begin{equation}
\label{original_milp}
\max_{(x_{prt}),\,(g_m,\tilde g_m)}
\quad
\sum_{p \in P}
\sum_{r \in R}
\sum_{t \in T}
k_p x_{prt}
-
\alpha
\sum_{m \in M}
\tilde g_m,
\end{equation}
where $\alpha > 0$ is a penalization parameter controlling the trade-off between patient prioritization and specialty-level balancing.

The complete \milp\ formulation is summarized below:
\begin{figure*}[t]
\begin{empheq}[left={{(\textit{MILP}-\atop \textit{PatAlloc})}\quad\empheqlbrace\quad}]{align*}    
\max \quad& \sum_{p\in P}\sum_{r\in R}\sum_{t\in T}k_p x_{prt} - \alpha \sum_{m\in M} \tilde{g}_m, \\[1mm]
s.t. \quad& \sum_{p \in P} x_{prt} (d_p + u) \leq c_{rt}, \quad \forall r \in R,\, t \in T \\[1mm]
\quad& \sum_{r\in R}\sum_{t\in T}x_{prt} \leq 1, \quad \forall p \in P, \\[1mm]
\quad& Y_m  = g_m + \frac{D_m}{D}\cdot \textit{Y}, \\[1mm]
\quad& \tilde{g}_m \geq g_m , \quad \tilde{g}_m \geq -g_m, \quad \forall m \in M. 
\end{empheq}
\end{figure*}

\section{Distribution of patient admissions}

To simulate the long-term evolution of the waiting list during the 20-week experiments, we analyzed the historical arrival process of new patients entering the elective surgery waiting list.

\begin{figure}[ht]
    \centering
    \includegraphics[width= \linewidth]{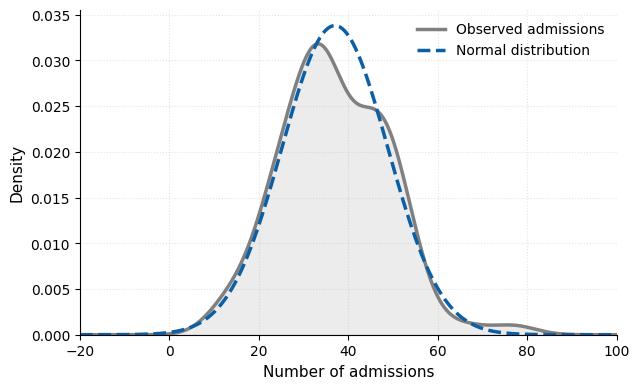}
    \caption{Empirical distribution of weekly patient admissions to the elective-surgery 
    waiting list together with the fitted Normal distribution.}
    \label{fig:admissionPatients}
\end{figure}

The number of weekly patient admissions was observed to follow an approximately Normal distribution, as illustrated in Figure~\ref{fig:admissionPatients}.

The fitted distribution has:
\begin{itemize}
    \item mean $\mu = 37$ weekly admissions;
    \item and standard deviation $\sigma = 11.8$.
\end{itemize}

During the simulations, the number of new patients entering the waiting list each week was generated according to this estimated distribution.

To preserve the heterogeneity observed in real hospital demand, new admissions were further distributed according to the empirical proportions of surgery types observed in the historical database. In total, $222$ distinct surgical procedures grouped into $10$ medical specialties were considered.

Table~\ref{tab:proporciones_operaciones} summarizes the empirical proportions associated with each surgical procedure and specialty. These proportions were used to generate the synthetic patient arrivals incorporated during the rolling-horizon simulations described in Section~\ref{sec:experimental_design}.

\begin{table*}
\centering
\resizebox{16cm}{!} {
\begin{tabular}{lrllrllrl}
\toprule
Operation &  Proportion &                        Service & Operation &  Proportion &                         Service & Operation &  Proportion &                         Service \\
\midrule
  1302029 &    0.231676 &  Otorhinolaryngological Surgery &      2104159 &      0.001593 &                         Traumatology &      1803011 &      0.000319 &                     Gastroenterology \\
   1902082 &    0.174634 &          Urology and Nephrology &      2104074 &      0.001593 &                         Traumatology &      1802002 &      0.000319 &                     Gastroenterology \\
   1902061 &    0.045889 &          Urology and Nephrology &      2104169 &      0.001593 &                         Traumatology &      1802036 &      0.000319 &                     Gastroenterology \\
   1802003 &    0.039834 &                Gastroenterology &      1802069 &      0.001593 &                     Gastroenterology &      1802068 &      0.000319 &                     Gastroenterology \\
   1302028 &    0.038560 &  Otorhinolaryngological Surgery &      1802024 &      0.001593 &                     Gastroenterology &      2104186 &      0.000319 &                         Traumatology \\
   1302008 &    0.028362 &  Otorhinolaryngological Surgery &      1302035 &      0.001593 &       Otorhinolaryngological Surgery &      1803015 &      0.000319 &                     Gastroenterology \\
   2703008 &    0.022626 &                       Dentistry &      1302061 &      0.001593 &       Otorhinolaryngological Surgery &      1802042 &      0.000319 &                     Gastroenterology \\
   1602202 &    0.016571 &       Dermatology and Teguments &      1703035 &      0.001275 &               Cardiovascular Surgery &      1502020 &      0.000319 &            Dermatology and Teguments \\
   2703004 &    0.015296 &                       Dentistry &      1302007 &      0.001275 &       Otorhinolaryngological Surgery &      1602001 &      0.000319 &            Dermatology and Teguments \\
   1502027 &    0.014659 &                 Plastic Surgery &      1802054 &      0.001275 &                     Gastroenterology &      1502035 &      0.000319 &                      Plastic Surgery \\
   1302012 &    0.012747 &  Otorhinolaryngological Surgery &      1902061 &      0.001275 &                         Traumatology &      1502056 &      0.000319 &                      Plastic Surgery \\
   1802014 &    0.010835 &                Gastroenterology &      1803031 &      0.001275 &  Proctolocic Rectal and Anal Surgery &      1502039 &      0.000319 &                      Plastic Surgery \\
   2106003 &    0.009879 &                    Traumatology &      1802032 &      0.001275 &                     Gastroenterology &      1502036 &      0.000319 &                      Plastic Surgery \\
   1902043 &    0.009560 &          Urology and Nephrology &      2104198 &      0.001275 &                         Traumatology &      1502032 &      0.000319 &                      Plastic Surgery \\
   1902013 &    0.009242 &          Urology and Nephrology &      2104101 &      0.001275 &                         Traumatology &      1402027 &      0.000319 &                Head and Neck Surgery \\
   1302023 &    0.008923 &  Otorhinolaryngological Surgery &      2104001 &      0.001275 &                         Traumatology &      1402042 &      0.000319 &                Head and Neck Surgery \\
   1602203 &    0.008923 &       Dermatology and Teguments &      1803016 &      0.001275 &                     Gastroenterology &      1402055 &      0.000319 &                Head and Neck Surgery \\
   1602231 &    0.007967 &       Dermatology and Teguments &      1803004 &      0.001275 &                     Gastroenterology &      1402002 &      0.000319 &                Head and Neck Surgery \\
   1502020 &    0.007648 &                 Plastic Surgery &      2104168 &      0.001275 &                         Traumatology &      1402051 &      0.000319 &                Head and Neck Surgery \\
   1502021 &    0.007648 &                 Plastic Surgery &      1502028 &      0.000956 &                      Plastic Surgery &      1402023 &      0.000319 &                Head and Neck Surgery \\
   1902064 &    0.007330 &          Urology and Nephrology &      1402017 &      0.000956 &                Head and Neck Surgery &      1402059 &      0.000319 &                Head and Neck Surgery \\
   2106001 &    0.007330 &                    Traumatology &      1302065 &      0.000956 &       Otorhinolaryngological Surgery &      1402038 &      0.000319 &                Head and Neck Surgery \\
   1402053 &    0.006692 &           Head and Neck Surgery &      1302042 &      0.000956 &       Otorhinolaryngological Surgery &      1402047 &      0.000319 &                Head and Neck Surgery \\
   1902082 &    0.006373 &                    Traumatology &      1802073 &      0.000956 &                     Gastroenterology &      1502050 &      0.000319 &                      Plastic Surgery \\
   1902022 &    0.006055 &          Urology and Nephrology &      2703023 &      0.000956 &                            Dentistry &      1502031 &      0.000319 &                      Plastic Surgery \\
   1302009 &    0.006055 &  Otorhinolaryngological Surgery &      1802031 &      0.000956 &                     Gastroenterology &      1902017 &      0.000319 &               Urology and Nephrology \\
   2104140 &    0.005736 &                    Traumatology &      1402006 &      0.000956 &                Head and Neck Surgery &      1902021 &      0.000319 &               Urology and Nephrology \\
   1902025 &    0.005736 &          Urology and Nephrology &      1803020 &      0.000956 &                     Gastroenterology &      1302022 &      0.000319 &       Otorhinolaryngological Surgery \\
   2104093 &    0.005417 &                    Traumatology &      1502014 &      0.000956 &                      Plastic Surgery &      1502020 &      0.000319 &       Otorhinolaryngological Surgery \\
   2104152 &    0.005417 &                    Traumatology &      1502063 &      0.000956 &                      Plastic Surgery &      1302059 &      0.000319 &       Otorhinolaryngological Surgery \\
   1802053 &    0.005099 &                Gastroenterology &      1902085 &      0.000956 &               Urology and Nephrology &      2104202 &      0.000319 &                         Traumatology \\
   1902060 &    0.005099 &          Urology and Nephrology &      1502003 &      0.000956 &                      Plastic Surgery &      2104201 &      0.000319 &                         Traumatology \\
   1302002 &    0.004780 &  Otorhinolaryngological Surgery &      1502016 &      0.000956 &                      Plastic Surgery &      2104014 &      0.000319 &                         Traumatology \\
   1302033 &    0.004780 &  Otorhinolaryngological Surgery &      2104107 &      0.000956 &                         Traumatology &      1902011 &      0.000319 &                         Traumatology \\
   1704004 &    0.004780 &          Cardiovascular Surgery &      2104003 &      0.000956 &                         Traumatology &      2104105 &      0.000319 &                         Traumatology \\
   2104057 &    0.004461 &                    Traumatology &      2104027 &      0.000956 &                         Traumatology &      2104141 &      0.000319 &                         Traumatology \\
   1402022 &    0.004143 &           Head and Neck Surgery &      1902016 &      0.000956 &               Urology and Nephrology &      1902032 &      0.000319 &                         Traumatology \\
   2104190 &    0.004143 &                    Traumatology &      1902018 &      0.000956 &               Urology and Nephrology &      2104073 &      0.000319 &                         Traumatology \\
   2104196 &    0.004143 &                    Traumatology &      1902090 &      0.000637 &               Urology and Nephrology &      1902064 &      0.000319 &                         Traumatology \\
   1502004 &    0.003824 &                 Plastic Surgery &      1402001 &      0.000637 &                Head and Neck Surgery &      2104199 &      0.000319 &                         Traumatology \\
   1502058 &    0.003824 &                 Plastic Surgery &      2104033 &      0.000637 &                         Traumatology &      1902073 &      0.000319 &                         Traumatology \\
   1802081 &    0.003824 &                Gastroenterology &      1402044 &      0.000637 &                Head and Neck Surgery &      1902044 &      0.000319 &                         Traumatology \\
   1402024 &    0.003505 &           Head and Neck Surgery &      1402052 &      0.000637 &                Head and Neck Surgery &      2104034 &      0.000319 &                         Traumatology \\
   1902075 &    0.003505 &          Urology and Nephrology &      2104111 &      0.000637 &                         Traumatology &      2104231 &      0.000319 &                         Traumatology \\
   1902044 &    0.003505 &          Urology and Nephrology &      2104197 &      0.000637 &                         Traumatology &      1902009 &      0.000319 &                         Traumatology \\
   1602224 &    0.003505 &       Dermatology and Teguments &      1302072 &      0.000637 &       Otorhinolaryngological Surgery &      2104108 &      0.000319 &                         Traumatology \\
   2104183 &    0.003187 &                    Traumatology &      2104026 &      0.000637 &                         Traumatology &      1902043 &      0.000319 &                         Traumatology \\
   1302003 &    0.003187 &  Otorhinolaryngological Surgery &      2104023 &      0.000637 &                         Traumatology &      1502039 &      0.000319 &       Otorhinolaryngological Surgery \\
   1802028 &    0.003187 &                Gastroenterology &      1703045 &      0.000637 &               Cardiovascular Surgery &      1302014 &      0.000319 &       Otorhinolaryngological Surgery \\
   1803031 &    0.003187 &                Gastroenterology &      1704040 &      0.000637 &               Cardiovascular Surgery &      1902040 &      0.000319 &               Urology and Nephrology \\
   1502025 &    0.002868 &                 Plastic Surgery &      1704050 &      0.000637 &               Cardiovascular Surgery &      1302071 &      0.000319 &       Otorhinolaryngological Surgery \\
   1302052 &    0.002868 &  Otorhinolaryngological Surgery &      1302063 &      0.000637 &       Otorhinolaryngological Surgery &      1902002 &      0.000319 &               Urology and Nephrology \\
   2104143 &    0.002868 &                    Traumatology &      1302049 &      0.000637 &       Otorhinolaryngological Surgery &      1902065 &      0.000319 &               Urology and Nephrology \\
   1802004 &    0.002549 &                Gastroenterology &      1802018 &      0.000637 &                     Gastroenterology &      1902019 &      0.000319 &               Urology and Nephrology \\
   1602002 &    0.002231 &       Dermatology and Teguments &      1502018 &      0.000637 &                      Plastic Surgery &      1902062 &      0.000319 &               Urology and Nephrology \\
   2104146 &    0.002231 &                    Traumatology &      2104151 &      0.000637 &                         Traumatology &      2104194 &      0.000319 &                         Traumatology \\
   1602223 &    0.002231 &       Dermatology and Teguments &      2104099 &      0.000637 &                         Traumatology &      1902047 &      0.000319 &               Urology and Nephrology \\
   2104136 &    0.002231 &                    Traumatology &      1502002 &      0.000637 &                      Plastic Surgery &      1902009 &      0.000319 &               Urology and Nephrology \\
   1302046 &    0.002231 &  Otorhinolaryngological Surgery &      1902067 &      0.000637 &               Urology and Nephrology &      1902011 &      0.000319 &               Urology and Nephrology \\
   1803022 &    0.002231 &                Gastroenterology &      1902048 &      0.000637 &               Urology and Nephrology &      1704003 &      0.000319 &               Cardiovascular Surgery \\
   1902029 &    0.001912 &          Urology and Nephrology &      1602201 &      0.000637 &            Dermatology and Teguments &      1704056 &      0.000319 &               Cardiovascular Surgery \\
   1704031 &    0.001912 &          Cardiovascular Surgery &      1602014 &      0.000637 &            Dermatology and Teguments &      1704005 &      0.000319 &               Cardiovascular Surgery \\
   1502047 &    0.001912 &                 Plastic Surgery &      1602206 &      0.000637 &            Dermatology and Teguments &      1704006 &      0.000319 &               Cardiovascular Surgery \\
   2104002 &    0.001912 &                    Traumatology &      2104189 &      0.000637 &                         Traumatology &      1402042 &      0.000319 &               Cardiovascular Surgery \\
   1902066 &    0.001912 &          Urology and Nephrology &      1802050 &      0.000637 &                     Gastroenterology &      1704064 &      0.000319 &               Cardiovascular Surgery \\
   1902073 &    0.001912 &          Urology and Nephrology &      1802058 &      0.000319 &                     Gastroenterology &      1704032 &      0.000319 &               Cardiovascular Surgery \\
   1902024 &    0.001912 &          Urology and Nephrology &      1802063 &      0.000319 &                     Gastroenterology &      1703021 &      0.000319 &               Cardiovascular Surgery \\
   2106002 &    0.001593 &                    Traumatology &      1802060 &      0.000319 &                     Gastroenterology &      1302057 &      0.000319 &       Otorhinolaryngological Surgery \\
   1302055 &    0.001593 &  Otorhinolaryngological Surgery &      1802101 &      0.000319 &                     Gastroenterology &      1302004 &      0.000319 &       Otorhinolaryngological Surgery \\
   1902050 &    0.001593 &          Urology and Nephrology &      1803019 &      0.000319 &                     Gastroenterology &      1502025 &      0.000319 &       Otorhinolaryngological Surgery \\
   1902084 &    0.001593 &          Urology and Nephrology &      1802080 &      0.000319 &                     Gastroenterology &      1803004 &      0.000319 &  Proctolocic Rectal and Anal Surgery \\
\bottomrule
\end{tabular}}
\caption{Proportion of weekly patient admissions by surgical procedure estimated from historical waiting-list records. }
\label{tab:proporciones_operaciones}
\end{table*}

\end{appendices}

\bibliography{ref}%

@ARTICLE{Zhou2020,
  author={Zhou, Liping and Geng, Na and Jiang, Zhibin and Wang, Xiuxian},
  journal={IEEE Transactions on Automation Science and Engineering}, 
  title={Public Hospital Inpatient Room Allocation and Patient Scheduling Considering Equity}, 
  year={2020},
  volume={17},
  number={3},
  pages={1124-1139},
  doi={10.1109/TASE.2019.2942990}}

@article{Roshanaei2020b,
title = {Reformulation, linearization, and decomposition techniques for balanced distributed operating room scheduling},
journal = {Omega},
volume = {93},
pages = {102043},
year = {2020},
issn = {0305-0483},
doi = {https://doi.org/10.1016/j.omega.2019.03.001},
author = {Vahid Roshanaei and Curtiss Luong and Dionne M. Aleman and David R. Urbach}
}

@article{Kamran2018,
title = {Uncertainty in advance scheduling problem in operating room planning},
journal = {Computers \&  Industrial Engineering},
volume = {126},
pages = {252-268},
year = {2018},
issn = {0360-8352},
doi = {https://doi.org/10.1016/j.cie.2018.09.030},
url = {https://www.sciencedirect.com/science/article/pii/S0360835218304467},
author = {Mehdi A. Kamran and Behrooz Karimi and Nico Dellaert},
}

@article{Azar2022,
title = {Dealing with uncertain surgery times in operating room scheduling},
journal = {European Journal of Operational Research},
volume = {299},
number = {1},
pages = {377-394},
year = {2022},
issn = {0377-2217},
doi = {https://doi.org/10.1016/j.ejor.2021.09.010},
author = {Macarena Azar and Rodrigo A. Carrasco and Susana Mondschein}
}

@article{Roshanaei2021,
author = {Roshanaei, Vahid and Naderi, Bahman},
doi = {10.1016/J.EJOR.2020.12.004},
issn = {0377-2217},
journal = {European Journal of Operational Research},
month = {aug},
number = {1},
pages = {65--78},
publisher = {North-Holland},
title = {{Solving integrated operating room planning and scheduling: Logic-based Benders decomposition versus Branch-Price-and-Cut}},
volume = {293},
year = {2021}
}

@article{Rizk2011,
author = {Rizk, Charbel and Arnaout, Jean-Paul},
doi = {10.1007/S10916-010-9648-Z},
issn = {1573-689X},
journal = {Journal of Medical Systems 2011 36:3},
month = {jan},
number = {3},
pages = {1891--1899},
publisher = {Springer},
title = {{ACO for the Surgical Cases Assignment Problem}},
url = {https://link.springer.com/article/10.1007/s10916-010-9648-z},
volume = {36},
year = {2011}
}

@article{Aringhieri2015,
author = {Aringhieri, Roberto and Landa, Paolo and Soriano, Patrick and T{\`{a}}nfani, Elena and Testi, Angela},
doi = {10.1016/J.COR.2014.08.014},
issn = {0305-0548},
journal = {Computers \&  Operations Research},
month = {feb},
pages = {21--34},
publisher = {Pergamon},
title = {{A two level metaheuristic for the operating room scheduling and assignment problem}},
volume = {54},
year = {2015}
}

@article{Rahimi2020,
author = {Rahimi, Iman and Gandomi, Amir H.},
doi = {10.1007/s11831-020-09432-2},
issn = {18861784},
journal = {Archives of Computational Methods in Engineering},
month = {may},
number = {3},
pages = {1667--1688},
publisher = {Springer},
title = {{A Comprehensive Review and Analysis of Operating Room and Surgery Scheduling}},
volume = {28},
year = {2020}
}

@article{Santibanez2007,
author = {Santib{\'{a}}{\~{n}}ez, Pablo and Begen, Mehmet and Atkins, Derek},
doi = {10.1007/s10729-007-9019-6},
issn = {13869620},
journal = {Health Care Management Science},
month = {sep},
number = {3},
pages = {269--282},
pmid = {17695137},
publisher = {Springer},
title = {{Surgical block scheduling in a system of hospitals: An application to resource and wait list management in a British Columbia health authority}},
volume = {10},
year = {2007}
}

@article{Jebali2006,
  title={Operating rooms scheduling},
  author={Jebali, Aida and Alouane, Atidel B Hadj and Ladet, Pierre},
  journal={International Journal of Production Economics},
  volume={99},
  number={1-2},
  pages={52--62},
  year={2006},
  doi = {10.1016/j.ijpe.2004.12.006},
  publisher={Elsevier}
}

@article{Atighehchian2020,
author = {Atighehchian, Arezoo and Sepehri, Mohammad Mehdi and Shadpour, Pejman and Kianfar, Kamran},
doi = {10.1007/s10479-019-03353-5},
issn = {15729338},
journal = {Annals of Operations Research},
month = {sep},
number = {1},
pages = {191--214},
publisher = {Springer},
title = {{A two-step stochastic approach for operating rooms scheduling in multi-resource environment}},
volume = {292},
year = {2020}
}

@article{Wang2020,
author = {Wang, Jin and Guo, Hainan and Tsui, Kwok Leung},
doi = {10.1080/00207543.2020.1815887},
issn = {1366588X},
journal = {International Journal of Production Research},
publisher = {Taylor and Francis Ltd.},
title = {{Two-stage robust optimisation for surgery scheduling considering surgeon collaboration}},
year = {2020}
}

@article{ThomasSchneider2020,
author = {{Thomas Schneider}, A. J. and {Theresia van Essen}, J. and Carlier, Mijke and Hans, Erwin W.},
doi = {10.1016/j.ejor.2019.09.029},
issn = {03772217},
journal = {European Journal of Operational Research},
month = {apr},
number = {2},
pages = {741--752},
publisher = {Elsevier B.V.},
title = {{Scheduling surgery groups considering multiple downstream resources}},
volume = {282},
year = {2020}
}

@article{Luo2019,
author = {Luo, Yan Yan and Wang, Bing},
doi = {10.1109/ACCESS.2019.2926780},
issn = {21693536},
journal = {IEEE Access},
pages = {102820--102831},
publisher = {Institute of Electrical and Electronics Engineers Inc.},
title = {{A New Method of Block Allocation Used in Two-Stage Operating Rooms Scheduling}},
volume = {7},
year = {2019}
}

@article{Latorre2016,
author = {Latorre-N{\'{u}}{\~{n}}ez, Guillermo and L{\"{u}}er-Villagra, Armin and Marianov, Vladimir and Obreque, Carlos and Ramis, Francisco and Neriz, Liliana},
doi = {10.1016/j.cie.2016.05.016},
issn = {03608352},
journal = {Computers and Industrial Engineering},
month = {jul},
pages = {248--257},
publisher = {Elsevier Ltd},
title = {{Scheduling operating rooms with consideration of all resources, post anesthesia beds and emergency surgeries}},
volume = {97},
year = {2016}
}

@article{wolff2012,
   author = {Guillermo Durán and Pablo A. Rey and Patricio Wolff},
   doi = {10.1007/S10479-016-2172-X},
   issn = {1572-9338},
   issue = {2},
   journal = {Annals of Operations Research 2016 258:2},
   month = {4},
   pages = {395-414},
   publisher = {Springer},
   title = {Solving the operating room scheduling problem with prioritized lists of patients},
   volume = {258},
   year = {2016},
}

@article{barrera2020,
  author = {Barrera, Javiera and Carrasco, Rodrigo A. and Mondschein, Susana and Canessa, Gianpiero and Rojas-Zalazar, David},
  doi = {10.1007/s10479-018-3008-7},
  issn = {15729338},
  journal = {Annals of Operations Research},
  month = {mar},
  number = {1-2},
  pages = {501--527},
  publisher = {Springer},
  title = {{Operating room scheduling under waiting time constraints: the Chilean GES plan}},
  volume = {286},
  year = {2020}
}

@article{bartek2018,
author = {Bartek, Matthew A. and Saxena, Rajeev C. and Solomon, Stuart and Fong, Christine T. and Behara, Lakshmana D. and Venigandla, Ravitheja and Velgapudi, Kalyani and Nair, Bala G. and Lang, John D.},
doi = {10.1016/j.jamcollsurg.2018.07.317},
issn = {10727515},
journal = {Journal of the American College of Surgeons},
month = {oct},
number = {4},
pages = {S149},
publisher = {Elsevier BV},
title = {{Improving Operating Room Efficiency: A Machine Learning Approach to Predict Case-Time Duration}},
volume = {227},
year = {2018}
}

@article{master2017,
  author = {Master, Neal and Zhou, Zhengyuan and Miller, Daniel and Scheinker, David and Bambos, Nicholas and Glynn, Peter},
  doi = {10.1007/s41060-017-0055-0},
  issn = {23644168},
  journal = {International Journal of Data Science and Analytics},
  month = {aug},
  number = {1},
  pages = {35--52},
  publisher = {Springer},
  title = {{Improving predictions of pediatric surgical durations with supervised learning}},
  volume = {4},
  year = {2017}
}

@article{edelman2017,
  author = {Edelman, Eric R. and van Kuijk, Sander M.J. and Hamaekers, Ankie E.W. and de Korte, Marcel J.M. and van Merode, Godefridus G. and Buhre, Wolfgang F.F.A.},
  doi = {10.3389/fmed.2017.00085},
  issn = {2296858X},
  journal = {Frontiers in Medicine},
  number = {JUN},
  pages = {85},
  publisher = {Frontiers Media S.A.},
  title = {{Improving the prediction of total surgical procedure time using linear regression modeling}},
  volume = {4},
  year = {2017}
}

@inproceedings{azar2017,
  address = {Seeon Abbey, Germany},
  author = {Azar, Macarena and Barrera, Javiera and Carrasco, Rodrigo A. and Mondschein, Susana},
  booktitle = {Proc. of the 13th Workshop on Models and Algorithms for Planning and Scheduling Problems},
  title = {Operating room scheduling with variable procedure times},
  year = {2017}
}

@article{samudra2016,
author = {Samudra, Michael and {Van Riet}, Carla and Demeulemeester, Erik and Cardoen, Brecht and Vansteenkiste, Nancy and Rademakers, Frank E.},
doi = {10.1007/s10951-016-0489-6},
issn = {10946136},
journal = {Journal of Scheduling},
month = {oct},
number = {5},
pages = {493--525},
publisher = {Springer New York LLC},
title = {{Scheduling operating rooms: achievements, challenges and pitfalls}},
volume = {19},
year = {2016}
}

@article{denton2007,
author = {Denton, Brian and Viapiano, James and Vogl, Andrea},
doi = {10.1007/s10729-006-9005-4},
issn = {13869620},
journal = {Health Care Management Science},
month = {feb},
number = {1},
pages = {13--24},
pmid = {17323652},
publisher = {Springer},
title = {{Optimization of surgery sequencing and scheduling decisions under uncertainty}},
url = {https://link.springer.com/article/10.1007/s10729-006-9005-4},
volume = {10},
year = {2007}
}

@article{dolatkhah2026reinforcement,
  title={A reinforcement-learning-based column generation algorithm for integrated operating room planning and scheduling},
  author={Dolatkhah, Mahdi and Hashemi Doulabi, Hossein and Rei, Walter and Gendreau, Michel},
  journal={International Journal of Production Research},
  pages={1--41},
  year={2026},
  publisher={Taylor \& Francis},
  doi={10.1080/00207543.2026.2637778}
}

@article{cardoen2010operating,
  title={Operating room planning and scheduling: A literature review},
  author={Cardoen, Brecht and Demeulemeester, Erik and Beli{\"e}n, Jeroen},
  journal={European journal of operational research},
  volume={201},
  number={3},
  pages={921--932},
  year={2010},
  publisher={Elsevier},
  doi={10.1016/j.ejor.2009.04.011}
}

@article{guerriero2011operational,
  title={Operational research in the management of the operating theatre: a survey},
  author={Guerriero, Francesca and Guido, Rosita},
  journal={Health care management science},
  volume={14},
  number={1},
  pages={89--114},
  year={2011},
  publisher={Springer},
  doi={10.1007/s10729-010-9143-6}
}

@article{zhu2019operating,
  title={Operating room planning and surgical case scheduling: a review of literature},
  author={Zhu, Shuwan and Fan, Wenjuan and Yang, Shanlin and Pei, Jun and Pardalos, Panos M},
  journal={Journal of Combinatorial Optimization},
  volume={37},
  number={3},
  pages={757--805},
  year={2019},
  publisher={Springer},
  doi={10.1007/s10878-018-0322-6}
}

@article{al2025comprehensive,
  title={A comprehensive review on operating room scheduling and optimization},
  author={Al Amin, Md and Baldacci, Roberto and Kayvanfar, Vahid},
  journal={Operational Research},
  volume={25},
  number={1},
  pages={3},
  year={2025},
  publisher={Springer},
  doi={10.1007/s12351-024-00884-z}
}
\end{document}